\input amstex
\documentstyle{amsppt}
\pageno=1
\magnification1200
\catcode`\@=11
\def\logo@{}
\catcode`\@=\active
\NoBlackBoxes
\vsize=23.5truecm
\hsize=16.5truecm

\def\slint{{\int\!\!\!\!\!\!@!@!@!@!{-}}}
\def\tslint{{\int\!\!\!\!\!@!@!@!{-}}}
\def\d{d\!@!@!@!@!@!{}^{@!@!\text{\rm--}}\!}
\def\supp{\operatorname{supp}}
\def\leg{\;\dot{\le}\;}

\def\crp{\overline{\Bbb R}_+}
\def\crm{\overline{\Bbb R}_-}

\def\rn{{\Bbb R}^n}
\def\rnp{{\Bbb R}^n_+}

\def\ang#1{\langle {#1} \rangle}

\def\rpp{ {\Bbb R}^2_{++}}
\def\rp{ \Bbb R_+}
\def\rmi{ \Bbb R_-}

\define\tr{\operatorname{tr}}
\define\Tr{\operatorname{Tr}}
\define\res{\operatorname{res}}
\define\TR{\operatorname{TR}}
\define\OP{\operatorname{OP}}
\define\deln{\tfrac{\partial  _\lambda ^{N-1}}{(N-1)!}}

\document

\topmatter
\title The local and global parts of the
basic zeta coefficient for operators on manifolds with boundary 
\endtitle
\author
Gerd Grubb
\endauthor
\affil
{Copenhagen Univ\. Math\. Dept\.,
Universitetsparken 5, DK-2100 Copenhagen, Denmark.
E-mail {\tt grubb\@math.ku.dk}}\endaffil
\subjclass 35S15, 58J42 \endsubjclass
\abstract
For operators on a compact manifold $X$ with boundary $\partial X$,
the basic zeta coefficient $C_0(B, P_{1,T})$ is the regular value at $s=0$ of the zeta function
$\Tr(B P_{1,T}^{-s})$, where $B=P_++G$ is a pseudodifferential
boundary operator (in the Boutet de Monvel calculus) --- for example
the solution operator of a classical elliptic problem --- and $P_{1,T}$ is
a realization of an elliptic differential operator $P_1$, having a ray
free of eigenvalues. 

Relative formulas (e.g.\ for the difference between the constants with
two different choices of $P_{1,T}$) have been known for some time and
are local. We here determine $C_0(B, P_{1,T})$ itself (when $P_1$ is of
even order), showing how it
is put together of local residue-type integrals (generalizing the
noncommutative residues of Wodzicki, Guillemin, Fedosov-Golse-Leichtnam-Schrohe) and
global canonical trace-type integrals (generalizing the canonical trace of Kontsevich
and Vishik, formed of Hadamard finite parts).

Our formula generalizes that of Paycha and Scott, shown recently for
manifolds without boundary. It leads in particular to new definitions
of noncommutative residues of expressions involving $\log P_{1,T}$.

Since the complex powers of $P_{1,T}$ lie far outside the Boutet de
Monvel calculus, the standard consideration of holomorphic families 
is not really useful here; instead we have developed a resolvent parametric 
method, where results from our calculus of parameter-dependent 
boundary operators can be used.

\endabstract
\rightheadtext{Basic zeta coefficient}
\endtopmatter

\subhead Introduction \endsubhead

The value of the zeta function at $s=0$ plays an important role in the
analysis of geometric invariants of operators on manifolds. For the
zeta function $\zeta (P_1,s)=\Tr P_1^{-s}$ (extended meromorphically to
${\Bbb C}$) defined from a classical elliptic pseudodifferential 
operator ($\psi
$do) $P_1$ on a closed manifold $X$, having a ray free of eigenvalues,
the value at $s=0$ is a fundamental ingredient in
index formulas. For the generalized zeta function  $\zeta
(A,P_1,s)=\Tr(AP_1^{-s})$, there is a pole at $s=0$ and the regular
value behind it serves as a ``regularized trace'' or ``weighted
trace'' 
(cf.\ e.g.\  Melrose et al.\ \cite{MN, MMS}, Paycha et al.\  \cite{CDMP, CDP}); it is likewise important in index
formulas.

Much is known for the case of closed manifolds:
The residue of $\zeta (A,P_1,s)$ at 0 
is proportional to Wodzicki's
noncommutative residue of $A$ (\cite{W}, see also Guillemin
\cite{Gu}). The regular value
at 0 (which we call {\it the basic zeta value}) equals
the Kontsevich-Vishik \cite{KV} canonical trace in special cases, and
in general there are defect formulas for it (formulas relating two different
choices of the auxiliary operator $P_1$, and formulas where $A$ is a
commutator), in terms of noncommutative
residues of related expressions involving $\log P_1$. The basic zeta
coefficient itself has recently been shown by Paycha and Scott
\cite{PS} to satisfy a formula 
with elements of canonical trace-type integrals (finite-part integrals
in the sense of Hadamard) defined from $A$, as well as
noncommutative residue-type
integrals defined from $A\log P_1$. 
The finite-part integral contributions are {\it global}, in the sense
that they depend on the full operator;
the residue-type contributions are {\it local}, in the sense
that they depend only on the strictly homogeneous symbols down to a
certain order.

For a compact manifold $X$ with boundary
$\partial X=X'$, one can investigate the analogous questions. $P_1$ is
here replaced by a realization $P_{1,T}$ of an elliptic differential
operator $P_1$ provided with a boundary condition $Tu=0$ at $X'$ such
that the resolvent $(P_{1,T}-\lambda )^{-1}$ exists on a ray in ${\Bbb
C}$, say ${\Bbb R}_-$. A suitable framework for these operators is the
Boutet de Monvel calculus (\cite{B}, \cite{G1}), which contains the
direct operators as well as their inverses and is closed under
composition of the various operator types (acting between bundles over
$X$ and $X'$). Now $A$ is replaced by an operator $B$ in the Boutet de
Monvel calculus, typically of the form $B=P_++G$, where $P$ is a $\psi
$do on a larger manifold $\widetilde X\supset X$, truncated to $X$,
and $G$ is a so-called singular Green operator (s.g.o.). For example,
$B$ can be the operator $\Delta _D^{-1}$ solving the Dirichlet problem
for a strongly elliptic differential operator $\Delta $.

The Boutet de Monvel calculus (the calculus of pseudodifferential
boundary operators, $\psi $dbo's) is a narrow calculus in the sense that
it has just what it takes to include elliptic differential boundary
problems and their solution operators in an ``algebra''. It requires
the $\psi $do component to satisfy the {\it transmission condition} at
the boundary (assuring that $C^\infty (X)$ is mapped into $C^\infty
(X)$), and is limited to $\psi $do's of integer order. There exist
other calculi for manifolds with boundary, such as e.g.\ the
b- or c-calculi of Melrose (\cite{M}, \cite{MN}), but these calculi
involve a degeneracy of the operators at the boundary, and do not
contain $\Delta _D^{-1}$. ---
One could say that the Boutet de Monvel calculus is more
``noncommutative'' than the other calculi, since it does {\it not} have the
property (typical of purely $\psi $do calculi) that
the commutator of two scalar operators $A_1$ and $A_2$ is of order $<
\operatorname{ord}A_1+\operatorname{ord}A_2$.

The celebrated method for determining spectral invariants such as
noncommutative residues and weighted traces, in the case of operators
on closed manifolds, or operators with a degeneracy at the boundary,
has been to study {\it holomorphic families}, i.e., families of operators depending holomorphically on
their order $z\in{\Bbb C}$ (\cite{Gu}, \cite{W}, \cite{KV}, \cite{MN},
\cite {CDMP}, \cite{PS}\dots). This does not work in the Boutet de Monvel
calculus, simply because the transmission condition in its
strict form allows integer orders only. Instead, we have had to
establish an alternative method: {\it the resolvent parametric method}.

The most frequently used holomorphic family is the family of powers;
in the present case it would be $BP_{1,T}^{-s}$. It can be defined
by Cauchy integrals from the resolvent parametric family 
$B(P_{1,T}-\lambda )^{-1}$, and the
invariants related to the holomorphic family $BP_{1,T}^{-s}$ match
certain invariants associated with the family
$B(P_{1,T}-\lambda )^{-1}$. In our method, we can bypass the detailed
analysis of the family $P_{1,T}^{-s}$ and carry out all the work on 
\linebreak$B(P_{1,T}-\lambda )^{-1}$ and its iterates $B(P_{1,T}-\lambda
)^{-N}$;
the only other operator we do have to 
involve is $\log P_{1,T}$. A detailed study of this operator
is worked out in \cite{GG}, joint with Gaarde. (In the holomorphic
calculi, the log-operator comes into the picture when
$z$-differentiation is performed; the reasons for its appearance in
the resolvent considerations are more subtle, connected with homogeneity
properties of symbols.) 

Previous results on noncommutative residues and zeta values for $\psi
$dbo's are as
follows: The noncommutative residue $\operatorname{res}(B)$ was 
introduced by Fedosov,
Golse, Leichtnam and Schrohe \cite{FGLS} as a sum of integrals of
symbol terms of order ``minus the dimension''. The
zeta function  $\zeta (B,P_{1,T},s)$ was defined
as the meromorphic extension of $\Tr(BP_{1,T}^{-s})$ to ${\Bbb C}$ 
by Grubb and Schrohe \cite{GSc1} (under restrictive hypotheses on
$P_{1,T}$), and the residue at $s=0$ identified with $\res B$.
``Weighted trace'' properties
(defect formulas) for the regular value at zero (with $P_{1,T}^{-s}$
replaced by its $\psi $do part) were 
established in \cite{GSc2,
G4}, including the fact that the value {\it modulo local terms} is a
finite-part integral defined from $B$. 

In the present paper we shall derive an explicit formula for
the regular value $C_0(B, P_{1,T})$ itself, with ingredients of the form of finite-part 
integrals (canonical trace-type terms) as well as residue-type
integrals involving $\log P_{1,T}$. The formula generalizes that of \cite{PS}, but has several more terms due
to the presence of the boundary.

The study leads to the introduction of a number of new noncommutative residue
formulas, generalizing those of \cite{W}, \cite{Gu} and \cite{FGLS}.
In particular, we show that the new residue definitions have a certain
traciality, vanishing on (suitable) commutators. 

  \subhead Plan of the paper \endsubhead Section 1 explains the
problem in more detail, desribing the asymptotic expansion of $\Tr
B(P_{1,T}-\lambda )^{-N}$, 
 and the decomposition of the relevant
coefficient in five terms (where $B=P_++G$ and $(P_{1,T}-\lambda )^{-1}=Q_{\lambda ,+}+G_\lambda $) 
$$
C_0(B, P_{1,T})=l_0((PQ_{\lambda})_+)- l_0(L(P,Q_\lambda ))+l_0(GQ_{\lambda,+} )+l_0(P_+G_\lambda )+l_0(GG_\lambda ).\tag0.1
$$
Section 2 establishes two general theorems describing the ``positive
regularity'' method initiated in \cite{G4} (based on \cite{G1}), 
and applies them to
complete the proofs of the needed asymptotic trace expansions. 
In Section 3, the nonlocal term $l_0((PQ_{\lambda })_+)$
is found by a method for closed manifolds originating from
\cite{G5}, and the local term 
$ l_0(L(P,Q_\lambda ))$ is determined by use of Section 2 combined
with a precise analysis. Section 4 treats the difficult nonlocal term
$l_0(GQ_{\lambda,+} )$, for which a new method involving Laguerre
expansions is introduced. Section 5 shows the residual nature of 
$l_0(P_+G_\lambda )$ and $l_0(GG_\lambda )$. In Section 6, we collect
the results in a theorem exhibiting the ingredients of canonical
trace-type and residue-type, extended from the case where
$P_1$ is of order 2 to general even-order cases. Here we moreover
establish some new defect formulas, relating two different choices of
$P_1$, or two different choices of branch cut of the logarithm, and
show the vanishing of the residues on certain commutators.

\subheading{1. Presentation of the problem and notation}

Consider a compact $n$-dimensional $C^\infty $ manifold $X$ with
boundary $\partial X=X'$, and a hermitian $C^\infty $ $M$-dimensional 
vector bundle
$E$ over $X$.
Let $B=P_++G$ be an operator of order $\sigma $ belonging to the
calculus of Boutet de Monvel \cite{B}, acting on sections of $E$. 
Here $P$ is a classical pseudodifferential operator, given on a larger
boundaryless $n$-dimensional manifold $\widetilde X$ in which $X$ is
smoothly imbedded, and acting in a bundle $\widetilde E$ extending $E$.
$P$ 
satisfies the transmission 
condition at $\partial X$, assuring that the truncation $P_+=r^+Pe^+$
preserves $C^\infty (X)$ (here $r^+$ restricts to $X^\circ$ and $e^+$
extends by zero on $\widetilde X\setminus X^\circ$). $G$ is a singular Green operator
(s.g.o.)\ of
class 0 with polyhomogeneous symbol. (More details on the calculus 
can be found e.g.\
in \cite{B} and Grubb \cite{G1}). When $P\ne 0$, we must assume
$\sigma \in\Bbb Z$ because of the requirements of the transmission
condition; when $P=0$, all  $\sigma \in\Bbb R$ are allowed. 

Along with $B$ we consider an auxiliary elliptic differential system
$\{P_1,T\}$ where $P_1$ is an elliptic differential operator of order
$m>0$ acting in $\widetilde E$, and $T$ is a differential trace operator, defining the 
realization $P_{1,T}$ in $L_2(X,E)$; its domain
$D(P_{1,T})$ consists of the
sections $u$ in the Sobolev space $H^m(X,E)$ with $Tu=0$. Here we assume that,
in local coordinates, the principal symbol $p^0_{1}(x,\xi
)$ has no eigenvalues on $\crm$, and the principal boundary symbol
operator 
$\{p^0_{1}(x',0,\xi
',D_n)- \lambda ,\, t^0(x',\xi ',D_n)\}$ 
is bijective for $\lambda \in
\crm$. Then $(p^0_{1}(x,\xi )-\lambda )^{-1}$ and the inverse of 
$\{p^0_{1}(x',0,\xi
',D_n)- \lambda ,\, t^0(x',\xi ',D_n)\}$ 
are defined for $\lambda $ in a
sector $V$ around 
$\rmi$, for all $x$, all $\xi '$ with $|\xi '|+|\lambda |\ne 0$.
It can be assumed that $\widetilde X$ is compact (cf.\ e.g.\ \cite{G2,
Pf.\ of Th.\ 7.4}).
The resolvent is$$
\aligned
R_\lambda &=(P_{1,T}-\lambda )^{-1}=Q_{\lambda ,+}+G_\lambda ,\text{
where }\\
Q_\lambda &=(P_1-\lambda )^{-1} \text{ on }\widetilde X.
\endaligned\tag 1.1
$$
These operators are defined except for $\lambda $ in a discrete subset of
$\Bbb C$; in particular they exist for large $\lambda $ in the sector $V$.
We can assume that no eigenvalues lie on $\rmi$ (by a rotation if
necessary), so that $\rmi$ is a so-called spectral cut.

The composed operator $BR_\lambda $ is trace-class when $m>\sigma +n$,
and we are interested in the expansion of its trace in powers of
$\lambda $ (with logarithmic factors).
If one does not want to assume that $P_1$ has a high order, one can instead
work with $N$'th powers of the resolvent.
Here we write
$$
\aligned
R^N_\lambda&=\deln R_\lambda =Q^N_{\lambda ,+}+G^{(N)}_\lambda ,\text{ with}\\
Q^N_{\lambda ,+}&=(Q^N_\lambda )_+=\deln Q_{\lambda,+} ,\quad G^{(N)}_\lambda =\deln G_\lambda .
\endaligned
\tag1.2
$$
 For certain calculations, it is
convenient to take $P_1$ of order 
$m=2$, considering  $\Tr(BR^N_{\lambda})$
for $N$ so large that $BR^N_{\lambda}$ is trace-class, namely
$N>(\sigma +n)/2$. 

The point of departure for the analysis of Laurent coefficients of
$\zeta (B,P_{1,T},s)$ is the
existence of trace expansions for $\lambda \to\infty $ on rays in $V$,
with $\delta >0$:
$$
\Tr(BR^N_{\lambda })
=\sum_{0\le j <n+\sigma }a^{(N)}_{j}\,(-\lambda )
^{\frac{n+\sigma -j}m-N}+(a^{(N)\prime}_0\log(-\lambda
)+a^{(N)\prime\prime}_0)(-\lambda )^{-N}+O(\lambda  ^{-N-\delta }),
\tag1.3
$$ 
when $N>(\sigma +n)/m$.
The expansion (1.3) was established by
Grubb and Schrohe in \cite{GSc1} when $m=2$, $P_1$ is principally scalar near
$X'$ and $T$ defines the
Dirichlet condition. (In fact, a full expansion with powers of
$-\lambda $ going to $-\infty $ was shown there, but we shall not need
the lower order terms in the present paper.) We shall show that it
holds for general second-order realizations $P_{1,T}$ with a spectral
cut (Corollary 2.7). For $m>n+\sigma $, the
expansion with $N=1$ was 
partially established in \cite{G4}, namely with $R_\lambda $ replaced by its
$\psi $do part $Q_{\lambda ,+}$:
$$
\Tr(BQ_{\lambda ,+})
=\sum_{0\le j <n+\sigma}c_{j}\,(-\lambda )
^{\frac{n+\sigma -j}m-1}+(c'_0\log(-\lambda
)+c''_0)(-\lambda )^{-1}+O(\lambda  ^{-1-\delta }),
\tag1.4
$$
 and we shall show the
supplementing result for the s.g.o.-part further below in Section
2 (Corollary 2.5). 
The expansion (1.4) (as well as expansions with $Q_{\lambda ,+}$ replaced by
$(Q^N_\lambda )_+$), 
has an interest in itself. 



The coefficients of $(-\lambda )^{-N}$ are particularly
interesting. In (1.3), one has for all $N>(\sigma
+n)/m$:
$$
a_0^{(N)\prime}=a'_0,\quad
a_0^{(N)\prime\prime}=a''_0-\alpha _Na'_0,\text{ with }
\alpha _N=\tsize\sum_{1\le j<N}\tfrac1j,
\tag1.5$$
where $a'_0$ and $a''_0$ are constants independent of $N$. 
A brief explanation of why
$\alpha _N$ enters is that it comes 
from derivatives of the log-term when $N$ is increased
(cf.\ (1.2)); 
a detailed account is given
in \cite{G5, Lemma 2.1}. It is the coefficient $a_0''$ that we are
searching for in particular.

The main efforts are made in the case where $P_1$ is of order 2; the
results are then extended to general cases of even $m$ in Section 6.

We do not here give a complete
treatment of cases where $m$ is odd.
Some results for such cases are included without extra effort (e.g.\
those in Section 5), for
other aspects we only give partial results; some indications of what may be
done in the first-order case are given in Remark 2.8 below. 
(Some difficulties are connected
with the fact that the classical part of $\log P_1$ does not satisfy
the transmission condition when $m$ is odd.)

To explain the connection with zeta functions, note that the
resolvents  can be used to define complex powers $P_1^{-s}$ 
 on $\widetilde X$ and $P_{1,T}^{-s}$ on $X$, by
Cauchy integral formulas such as
$$
P_{1,T}^{-s}=\tfrac i{2\pi }\int_{\Cal C}\lambda
^{-s} \,(P_{1,T}-\lambda )^{-1}\, d \lambda ,
\tag1.6
$$ 
where $\Cal C$ is a curve in ${\Bbb C}\setminus\crm$ encircling the
nonzero spectrum of $P_{1,T}$ in the positive direction (e.g.\ a Laurent loop around $\crm$). By
the transition formulas accounted for e.g.\ in \cite{GS2},
the above resolvent trace expansions imply that the zeta functions
$$
\zeta (B,P_{1,T},s)=\Tr(BP_{1,T}^{-s}),\quad \zeta (B,P_{1,+},s)
=\Tr(B(P_{1}^{-s})_+),\tag1.7
$$
holomorphic for $\operatorname{Re}s>(\sigma +n)/m$, extend
meromorphically to $\operatorname{Re}s>-\delta $, with simple poles at
the real points $(\sigma +n-j)/m$ and in particular a Laurent
expansion at $s=0$:
$$
\aligned
\zeta (B,P_{1,T},s)&= C_{-1}(B,P_{1,T})s^{-1}+[C_0(B,P_{1,T})-\Tr(B\Pi
_0(P_{1,T}))]s^0+O(s),\\
\zeta (B,P_{1,+},s)&=C_{-1}(B,P_{1,+})s^{-1}+[C_0(B,P_{1,+})-\Tr(B\Pi
_0(P_1)_+)]s^0+O(s),
\endaligned\tag1.8
$$
where $\Pi _0(P_1)$ resp.\ $\Pi _0(P_{1,T})$ is the generalized
eigenprojection for the zero eigenvalue of $P_1$ resp.\ $P_{1,T}$.
Here, in relation to (1.3)--(1.4) with (1.5),
$$
\aligned
 C_{-1}(B,P_{1,T})&=a_0'=a_0^{(N)\prime},\\
C_{-1}(B,P_{1,+})&=c_0',
\endaligned\tag1.9
$$
and
$$
\aligned
 C_{0}(B,P_{1,T})&=a_0''=a_0^{(N)\prime\prime}+\alpha _Na_0',\\
C_{0}(B,P_{1,+})&=c_0'',
\endaligned\tag1.10
$$
 The $\Pi _0$-terms enter only in the
zeta-expansions (1.8), not in the resolvent expansions
(1.3)--(1.4); they stem from the fact that 0 is excluded from the
contour integration (1.6).

It is known from \cite{GSc1} for (1.3) with $m=2$, from \cite{G4} for
(1.4), and will follow easily for general cases from the proofs given below, that
$$
C_{-1}(B, P_{1,T})=C_{-1}(B, P_{1,+})=\tfrac1m \operatorname{res}B,\tag1.11
$$
where  $\res B=\res(P_++G)$ is the noncommutative residue; the linear 
functional
defined by Fedosov, Golse, Leichtnam and Schrohe \cite{FGLS}:
$$
\res(P_+)=\int_X\int_{|\xi |=1}\tr p_{-n}(x,\xi )\,\d S(\xi ),\quad
\res G=\res_{X'}(\tr_n G)
$$
(expressed in local coordinates). Here $\res_{X'}$ is the Wodzicki
residue of the $\psi $do $\tr_nG$ over $X'$ (the normal trace of $G$, see
(2.2) below). The indication $\tr$ means fiber trace, $\d \xi $ stands
for $(2\pi )^{-n}\,
d\xi $, and $\d S(\xi )$ is $(2\pi )^{-n}$ times the usual surface measure.
 This 
coefficient is completely independent of $P_1$ or $T$, and is
local. It is zero when $\sigma +n\notin {\Bbb N}=\{0,1,2,\dots\}$.

We shall focus the attention on the coefficient $C_{0}(B,P_{1,T})$,
called {\it the basic zeta coefficient}, as well as its variant
$C_{0}(B,P_{1,+})$, and we search for explicit formulas for these constants.

In the papers \cite{GSc2}, \cite{G4} we have studied the
trace defect $C_0(B,P_{1,+})-C_0(B,P_{2,+})$, which compares the
coefficients for two different auxiliary operators $P_1$ and $P_2$.
It was shown to be local in
\cite{GSc2}. Moreover, \cite{G4} showed
that it can be expressed as a noncommutative residue in the sense of \cite{FGLS}:$$
C_0(B,P_{1,+})-C_0(B,P_{2,+})=-\tfrac 1m \operatorname{res}(B(\log
P_1-\log P_2)_+),\tag 1.12
$$ 
when $m$ is even. In \cite{G4}, $m>\sigma +n$ or $m=2$; general 
values of $m$
are covered by Theorem 2.6 below.
 When $m$ is odd, the residue definition of
\cite{FGLS} is not directly applicable, but there is a more
complicated residue-like 
interpretation of the right-hand side in (1.12).  

The operator $\log P_1$ can be 
defined (on smooth functions) by a Cauchy integral and approximation:
$$
\log P_1=\lim _{s\searrow 0}\tfrac i{2\pi }\int_{\Cal C}\lambda
^{-s}\log \lambda \,(P_1-\lambda )^{-1}\, d \lambda ,\tag1.13
$$ 
where $\Cal C$ is a curve in ${\Bbb C}\setminus\crm$ encircling the
nonzero spectrum of $P_1$. Its symbol in local coordinates is of the form
$$
\gathered
\operatorname{symb}(\log P_1)= m\log[\xi ]+l(x,\xi ),\text{ with $l(x,\xi )$ classical of order 0,}\\
[\xi ]\text{ is a  $C^\infty $ function $\ge \tfrac12 $ with }[\xi ]
=|\xi |\text{ for
}|\xi |\ge 1 ,
\endgathered
\tag1.14
$$ 
cf.\ e.g.\ Okikiolu \cite{O}, and $l(x,\xi )$ satisfies the
transmission condition when $m$ is even (details in
\cite{GG, Lemma 2.1}).
 
In the present paper, we shall prove a formula for $C_0(B,P_{1,+})$
itself, showing how it is, in 
local coordinates, the sum of a finite-part integral (in the sense of
Hadamard)
defined from $B$ and a residue-type
integral defined from 
$B(\log P_1)_+$, suitably interpreted. The value {\it modulo local
terms} was found in \cite{GSc2}; the present description is
much more precise. Moreover, we shall extend the description to 
$C_0(B,P_{1,T})$.

\medskip

Our result is a generalization of, and is inspired from, the recent result of
Paycha and Scott \cite{PS} on Laurent coefficients in the
boundaryless case, further analyzed in our note \cite{G5}. They showed the
following:

Let $A$ and $P_1$ be classical
pseudodifferential 
operators of order $\sigma \in\Bbb R$ resp.\ $m\in\Bbb R_+$ on a
compact $n$-dimensional manifold $\widetilde X$ without boundary, $P_1$ being
elliptic with $\Bbb R_-$ 
as a spectral cut. Then
all the coefficients in the Laurent expansions of the generalized
zeta function $\zeta (A,P_1,s)=\Tr(AP_1^{-s})$ around the poles can be
expressed as combinations 
of finite-part integrals and residue-type integrals of associated
logarithmic symbols. Denote (similarly to (1.8)) the Laurent coefficient of $s^0$ at
zero (the regular value at zero) by $C_0(A,P_1)-\Tr(A\Pi _0(P_1))$, where
$\Pi _0(P_1)$ is the 
generalized zero eigenprojection of $P_1$.
Then 
$C_0(A,P_1)$ satisfies
$$
C_0(A,P_1)=\int_{\widetilde X}\bigl(\TR_x(A)-\tfrac1m\operatorname{res}_{x,0}(A\log
P_1)\bigr)\,dx.\tag1.15
$$
Here the function $\TR_x(A)-\tfrac1m\operatorname{res}_{x,0}(A\log
P_1)$ is defined in a local coordinate system by:$$
\gather
\TR_x(A)=\slint
\tr a(x,\xi )\,\d\xi,\tag1.16\\
 \operatorname{res}_{x,0}(A\log P_1)
=\int_{|\xi |=1}\tr r_{-n,0}(x,\xi )\d S(\xi ).\tag1.17
\endgather
$$

For (1.16), the expression $\tslint
\tr a(x,\xi )\,\d\xi$ 
is defined for each $x$ 
as a Hadamard finite-part integral,
namely as the constant term in the asymptotic expansion 
of $\int_{|\xi |\le R }\tr a(x,\xi  )\,\d\xi $ in 
powers $R ^{\sigma +n-j}$ ($j\in\Bbb N$), $R ^0$ and $\log R $
(cf.\ Lesch \cite{L}, see also \cite{GSc2, (3.12)ff.}). The notation
$\TR_x$ is inspired from the notation of \cite{KV}; in fact, as
pointed out in \cite{L}, $\TR_xA$ 
integrates in suitable cases to the canonical trace $\TR A$
 (see also
\cite{G3}).

For (1.17), the symbol 
of $R=A\log P_1$ is denoted $r(x,\xi )$; it is log-polyhomogeneous of the form (cf.\ (1.14))$$
r(x,\xi 
)\sim \sum_{j\in\Bbb N}(r_{\sigma -j,0}(x,\xi )+r_{\sigma -j,1}(x,\xi )
\log[\xi ]),
\tag1.18$$
where each $r_{\sigma -j,l}$ is homogeneous in $\xi $ of degree $\sigma -j$
for $|\xi |\ge 1$. So $r_{-n,0}$ is the log-free term of order $-n$
(taken equal to 0 when $\sigma +n\notin\Bbb N$).

We find it interesting that the two ``competing'' trace functionals res
and TR both enter here, in localized versions (only the
collected expression
$\bigl(\TR_x(A)-\frac1m\operatorname{res}_{x,0}(A\log
P_1)\bigr)\,dx$ has an invariant meaning as a 
density on $\widetilde X$).
Note that $\res_{x,0}(A\log P_1)$ is {\it local} (depends on homogeneous
symbols down to order $-n$),
whereas $\TR_x(A)$ is {\it global} (depends on the full structure).

The formula was shown in \cite{PS} by use of holomorphic families of
$\psi $do's (depending holomorphically on their complex order $z$).
We showed in \cite{G5} how the formula could be derived by
methods relying directly on the knowledge of the resolvent
$(P_1-\lambda )^{-1}$, as a preparation for the present
generalization to manifolds with boundary.

Like \cite{PS}, we shall use the notation$$
\operatorname{res}_x(Q)=\int_{|\xi |=1}\tr q_{-n}(x,\xi )\,\d S(\xi ),
\tag1.19$$
when $Q$ is a classical $\psi $do on $\widetilde X$ with symbol
$q(x,\xi )$ in local coordinates. Here 
$\operatorname{res}_x(Q)$ has a meaning only in local 
coordinates, but the integral of
$\operatorname{res}_x(Q)$ over 
$\widetilde X$ can be given an invariant meaning 
as the Wodzicki noncommutative residue $\operatorname{res}Q$.

\example{Remark 1.1}
In view of (1.14), we can in local coordinates
 define $(\log P_1)_{\operatorname{cl}}=\OP(l(x,\xi ))$ as
``the classical part of $\log P_1$''. Then, as is easily checked from
the composition rules,
 $$
\operatorname{res}_{x,0}(A\log P_1)=\operatorname{res}_x(A(\log
P_1)_{\operatorname{cl}}),\tag1.20$$ 
when we use the notation  (1.19).
But this generally has a meaning only in local coordinates.
 \endexample

Our goal is to generalize (1.15) to a characterization of $C_0(B,
P_{1,T})$, via a study of the coefficient of $(-\lambda )^{-N}$ in
the trace expansion of $BR^N_ \lambda =\deln BR_\lambda $. The
operator breaks
up in five terms:
$$
\align
BR_\lambda &=(P_++G)(Q_{\lambda ,+} +G_\lambda )\tag1.21
\\
&=P_+Q_{\lambda,+}+GQ_{\lambda ,+}+(P_++G)G_\lambda ,
\\
&=(PQ_\lambda
)_+-L(P,Q_\lambda )+GQ_{\lambda ,+}+P_+G_\lambda +GG_\lambda ,
\\
BR^N_\lambda &=(PQ^N_\lambda
)_+-L(P,Q^N_\lambda )+GQ^N_{\lambda ,+}+P_+G^{(N)}_\lambda +GG^{(N)}_\lambda ,
\tag1.22\endalign
$$
where we use the notation (1.2), and the notation $L(P,P')$ for the
``leftover term'' in the truncation of a product:
$L(P,P')=(PP')_+-P_+P'_+$ (more on this in (3.16)ff.).
The part $BQ^N_{\lambda ,+}$ (the first three terms in (1.22))
has 
been studied
in \cite{GSc1}, \cite{GSc2} and \cite{G4}, the last paper giving 
 the trace defect
formula (1.12). For
this part, the strategy will be:
\roster

\item {\it Find $C_0(B, P_{1,+})$ for one particularly manageable choice of $P_1$.}

\item  {\it Extend to more general $P_2$ by combination with the
trace defect formula 
{\rm (1.12)}.}
\endroster
This program is carried out in Sections  3 and 4.
 
The part $BG^{(N)}_\lambda $ (the last two terms in (1.22)) has been
treated under special assumptions and with only qualitative results in
\cite{GSc1} and \cite{GSc2}. An exact trace defect formula was not
shown, so we would need to find this. In fact, we shall aim directly for a
precise description of the coefficient in question (from which a
defect formula follows), since this can be done with the
methods established in \cite{G4}. The coefficient is {\it
local} --- as shown under restrictive hypotheses in 
\cite{GSc1}.  Part of the study is carried out in Section 2, showing
the existence of an expansion with a qualitative description of the
coefficients, and the connection with our main problem is worked out
in Section 5.

It will be practical to introduce a general notation for the relevant
coefficient in various trace expansions (cf.\ (1.5)):

\proclaim{Definition 1.2} Let $M_{\lambda }$ be an operator family
defined for large $\lambda $ in a sector $V$ of ${\Bbb C}$, with
$M_{\lambda ,N}=\deln M_\lambda 
$ defined there for all positive integers
$N$.
If the $M_{\lambda ,N}$ for all $N>N _0$ have trace expansions of the form
$$\multline
\Tr M_{\lambda ,N}=\sum_{0\le j <n+\sigma }a_{N,j}\,(-\lambda )
^{\frac{n+\sigma -j}2-N}\\
+(a'_0\log(-\lambda
)+a''_0-\alpha _Na'_0)(-\lambda )^{-N}+O(\lambda  ^{-N-\delta
}),\text{ some }\delta >0,
\endmultline\tag1.23
$$
 for $\lambda \to\infty $ on rays in $V$, with $a'_0$ and $a''_0$
independent of $N$, $\alpha _N=\tsize\sum_{1\le j<N}\tfrac1j$, we define
$$
l_0(M_{\lambda })=a''_0.\tag 1.24
$$

\endproclaim

The main result is the following theorem with $P_1$ of even order $m$, 
whose ingredients will
be explained in the next sections:

\proclaim{Theorem 1.3}
The basic zeta coefficient $C_0(B, P_{1,T})$ is a sum of terms,
calculated in local coordinates:
$$
\aligned
C_0(B, P_{1,T})&=\int_X[\TR _xP-\tfrac1m\operatorname{res}_{x,0}(P\log
P_1)]\, dx+\tfrac1m\operatorname{res} L(P,\log P_1)
\\
&\quad+\int _{X'}[\TR _{x'}\tr_n
G-\tfrac1m\operatorname{res}_{x',0}\tr'_n(G(\log P_1)_+)]\, dx'\\ 
&\quad
-\tfrac1m\res(P_+G_1^{\log})-\tfrac1m\res(GG_1^{\log}).
\endaligned
\tag1.25$$
Here $C_0(B, P_{1,+})$ is the sum of the first three terms.
\endproclaim

The five terms are found as the coefficient $l_0$ in trace expansions of the
five terms in (1.22), as written in (0.1). The residues appearing in (1.25) are various 
generalizations of the definition of \cite{FGLS}, and will be suitably 
explained in the process of deduction of the formula. Both 
$\operatorname{res}L(P,\log P_1)$ and $\res(GG_1^{\log})$ are
Wodzicki-type residues over $X'$ of $\psi $do's $\tr_n L(P,\log
P_1)$ resp.  $\tr_n(GG_1^{\log})$, whereas  $\res(P_+G_1^{\log})$,
as well as $\operatorname{res}_{x',0}\tr'_n(G(\log P_1)_+)$ are more
delicate to define. The theorem is shown in Section 6, based on the
results from Sections 3 and 4 (with $m=2$) and 5 (with $m>0$); moreover,
new defect formulas  are derived. In
particular, it is proved that the new residue definition is
tracial in the sense that
$$
\res ([B, \log P_{1,T}])=0,\tag1.26
$$
when $B$ is of order and class 0.

\subhead 2. 
A strategy for local contributions
\endsubhead

In this section we show the supplementary results needed to get the
trace expansions (1.3) and (1.4) in the stated generality.
The terms $L(P,Q_\lambda ^N)$, 
$P_+G_\lambda ^{(N)}$ and $GG_\lambda ^{(N)}$
 in (1.22) will be seen to have trace
expansions (1.23) with $a'_0=0$, and $a''_0$ {\it local}. To show
this result, we can apply a method
introduced in \cite{G4}, relying on the ``positive regularity'' (in
the sense of \cite{G1}) of the $\psi $do family on $X'$ obtained by taking
normal traces. 

We shall explain the method here, establishing the main technical
points in a form applicable to general situations, and give its
application to $P_+G_\lambda ^{(N)}+GG_\lambda ^{(N)}=BG_\lambda
^{(N)}$ in the present section, treating $L(P,Q_\lambda ^N)$ later in Section
3.

The trace expansions are derived as finite sums of trace expansions
worked out in local coordinates, where the situation is carried over to
$\rnp$, with coordinates $x'\in {\Bbb R}^{n-1}$, $x_n\in \rp$. We
recall that when $G_0$ is a
singular Green operator on $\rnp$ of order $<1-n$ and class 0 with a kernel having
compact $(x',y')$-support (hence is trace-class), then $\Tr
G_0=\Tr_{{\Bbb R}^{n-1}}(\tr_n G_0)$, where $\tr_n$ indicates the normal
trace.
To explain this further, recall from \cite{G1} that when $G_0$ 
has the symbol $g_0(x',\xi ',\xi _n,\eta _n)$ 
and symbol-kernel $\tilde g_0(x',x_n,y_n,\xi ')=r^+_{x_n}r^+_{y_n}\Cal
F^{-1}_{\xi _n\to x_n}\overline{\Cal F}^{-1}_{\eta  _n\to
y_n}g_0(x',\xi ',\xi _n,\eta _n)$, the action of $G_0$ is defined for
$u\in r^+\Cal S({\Bbb R}^n)$ by 
$$
G_0u=\int _{\Bbb R^{n-1}}e^{ix'\cdot \xi '}\int_0^\infty
\tilde g_0(x',x_n,y_n,\xi ')\acute u(\xi ',y_n)\,dy_n\d \xi ',
\; \acute u=\Cal F_{x'\to \xi '}u
\tag2.1
$$
($r^+$ restricts to $\{ x_n>0\}$), and the normal trace $\tr_nG_0$ is a $\psi $do $S_0$ on $\Bbb R^{n-1}$ with symbol
$$
\aligned
s_0(x',\xi ')&=(\tr_ng_0)(x',\xi ')=\int_0^\infty \tilde g_0(x',x_n,x_n,\xi
' ) \, dx_n\\
&=\int_0^\infty r^+_{x_n}r^+_{y_n}\int^+\int ^+e^{ix_n\xi _n-ix_n\eta
_n}g_0(x',\xi ',\xi _n,\eta _n )\,\d\xi _n\d\eta _n dx_n
\\
&=\int^+ g_0(x',\xi ',\xi _n,\xi _n )\,\d\xi _n.
\endaligned
\tag2.2
$$
The plus-integral stands for an extension of the usual integral (cf.\
e.g.\ \cite{G1, p.\ 166}), and
the last formula follows since the factor $e^{ix_n\xi _n-ix_n\eta
_n}$ together with the integrations in $\eta _n$ and $x_n$ give rise
to a backwards and a forwards Fourier transform.

We denote
$
(1+|\xi '|^2)^{\frac12}=\ang{\xi '}$, $(1+|\xi
'|^2+\mu ^2)^{\frac12}=\ang{\xi ',\mu },
$  
and use the sign $\leg$ in inequalities to indicate ``$\le$ a
constant times''.

The idea is to show that taking the normal
trace of $BG_\lambda ^{(N)}$, in local coordinates,  leads to a $\lambda $-dependent family
of $\psi $do's on ${\Bbb R}^{n-1}$, where a 
kernel expansion and a formula for the basic
coefficient can be deduced exactly as in the proof of \cite{G4, Th.\
4.5}. For the sake of general applicability,
the crucial steps in that proof will be formulated in the 
following two general theorems, where we outline the proofs with
reference to \cite{G4}. 

\proclaim{Theorem 2.1} 
Let $\sigma \in \Bbb R$,
let $\delta \in \,]0,1[\,$ such that $\sigma +\delta \notin{\Bbb Z}$, 
let $m$ and $N$ be positive integers, and let $\Cal
S^{(N)}_\lambda =\OP(\frak s^{(N)}(x',\xi ',\lambda ))$ be a  
family of $\psi $do's on $\Bbb R^{n-1}$ depending holomorphically on
$\lambda $ in a keyhole region $V'=V\cup \{0<|\lambda |<r\}$,
$V$ being a sector around ${\Bbb R}_-$, of order $\sigma -mN$ 
and regularity
$\sigma -\delta $ in terms of $\mu $ on each ray $\lambda =-\mu
^me^{i\theta }$ in $V$. Assume moreover that $(-\lambda )^{N}\Cal
S^{(N)}_\lambda $ is of order $\sigma $ and 
regularity $\sigma +\delta $, and that the symbol satisfies
$$
\aligned
|\partial _{x',\xi '}^{\beta ,\alpha } [\frak s^{(N)}(x',\xi ',\lambda
)-\sum_{j<J}\frak s^{(N)}_{\sigma   -mN-j}(x',\xi ',\lambda )]|
&\leg \ang{\xi '}^{\sigma   - \delta  -|\alpha |-J}\ang{\xi ',\mu
}^{-mN+\delta },\\
|\partial _{x',\xi '}^{\beta ,\alpha } [\frak s^{(N)}(x',\xi ',\lambda
)-\sum_{j<J}\frak s^{(N)}_{\sigma   -mN-j}(x',\xi ',\lambda )]|
&\leg \ang{\xi '}^{\sigma   +\delta -|\alpha |-J}\ang{\xi ',\mu
}^{-\delta }\mu ^{-mN},
\endaligned
\tag2.3
$$
on the rays in $V$, for all $\alpha 
,\beta ,J$. 

 If $N>(\sigma +n-1)/m$,  the kernel on the diagonal has an
expansion
$$
\aligned
K(\Cal S^{(N)}_\lambda ,x',x')&=\sum_{0\le l\le \sigma +n-1}\tilde
{\frak s}^{(N)}_l(x')(-\lambda )
^{\frac{n-1+\sigma -l}m -N}+O(\lambda ^{-N-\frac\delta m}),\text{ with}\\
\tilde {\frak s}^{(N)}_{l}(x')&=
\int_{\Bbb R^{n-1}}{\frak s}^{(N)h}_{ \sigma -mN-l}(x',\xi ',-1)\,\d\xi',
\endaligned
\tag2.4
$$
where the strictly homogeneous symbols  
$\frak s^{(N)h}_{\sigma
-mN-l}$ are integrable at $\xi '=0$ for $l\le \sigma +n-1 $. If
$\sigma +n-1\in{\Bbb N}$, the coefficient of $(-\lambda ) ^{-N}$ is$$
\tilde {\frak s}^{(N)}_{\sigma +n-1 }(x')=
\int_{\Bbb R^{n-1}}{\frak s}^{(N)h}_{ -mN-n+1}(x',\xi ',-1)\,\d\xi'.\tag2.5
$$
If $\sigma +n-1\notin{\Bbb N}$, there is no term with $(-\lambda )
^{-N}$; we include one trivially by setting $\tilde {\frak s}^{(N)}_{\sigma +n-1 }(x')=
0$.\endproclaim

\demo{Proof} The strictly homogeneous version $\frak s^{(N)h}_{j}$ of
$\frak s^{(N)}_{j}$ is the extension of $\frak s^{(N)}_{j}|_{|\xi
'|\ge 1}$ by homogeneity into
$\{1>|\xi '|>0\}$.

The expansion down to the term with $l=n-2+\sigma $ is
assured by the general theory of \cite{G1}, cf.\ e.g.\ Lemmas 3.1--5
in \cite{G4}. To include the next term, one proceeds as in the proof
of \cite{G4, (4.35)ff.}, if $\sigma +n-1\in{\Bbb N}$ (with the number
$\frac14$  
replaced by $\delta >0$): 

The
estimates (2.3)
imply that the symbol with $l=\sigma +n-1 $ satisfies
$$\aligned
|\frak s^{(N)}_{-mN-n+1}(x',\xi ',\lambda )|&\leg \ang{\xi '}^{-\delta -n+1}
\ang{\xi ',\mu
}^{-mN+\delta },\\
|\frak s^{(N)}_{-mN-n+1}(x',\xi ',\lambda )|&\leg \ang{\xi '}^{\delta -n+1}
\ang{\xi ',\mu
}^{-\delta }\mu ^{-mN},
\endaligned\tag2.6$$
and the remainder $\frak s^{(N)\prime}=\frak s^{(N)}-\sum_{l< \sigma +n}\frak
s^{(N)}_{\sigma -mN-l+1}$ after this term satisfies
$$\aligned
|\frak s^{(N)\prime}|&\leg \ang{\xi '}^{-\delta -n}\ang{\xi ',\mu
}^{-mN+\delta },\\
|\frak s^{(N)\prime}|
&\leg \ang{\xi '}^{\delta -n}\ang{\xi ',\mu
}^{-\delta }\mu ^{-mN}.
\endaligned\tag2.6$$
From (2.6) follows as in \cite{G1, Lemma 2.1.9} that similar
estimates are valid for the strictly homogeneous symbols:
$$\aligned
|\frak s^{(N)h}_{-mN-n+1}(x',\xi ',\lambda )|&\leg |\xi '|^{-\delta -n+1}
|(\xi ',\mu)|^{-mN+\delta },\\
|\frak s^{(N)h}_{-mN-n+1}(x',\xi ',\lambda )|&\leg |\xi '|^{\delta -n+1}
|(\xi ',\mu )
|^{-\delta }\mu ^{-mN},
\endaligned\tag2.8$$
so $\frak s^{(N)h}_{-mN-n+1}$ is integrable at $\xi '=0$ (besides being so
for $|\xi '|\to \infty $) when $\lambda  \ne 0$. Then 
$$
\aligned
K(\operatorname{OP}(\frak s^{(N)h}_{-mN-n+1}),x',x')&=
\tilde {\frak s}^{(N)}_{\sigma +n-1}(x')\,(-\lambda )^{-1},\text{ with}\\
\tilde {\frak s}^{(N)}_{\sigma +n-1}(x')&=
\int_{\Bbb R^{n-1}}{\frak s}^{(N)h}_{-mN-n+1}(x',\xi ',-1)\,\d\xi',
\endaligned\tag2.9
$$
as desired. This gives the needed extra term. The remainder estimate
is obtained (like the preceding
lines) exactly as in \cite{G4, Th.\ 4.5}; we shall not repeat the
details. If $\sigma +n-1\notin{\Bbb N}$, there is no term with
$(-\lambda )^{-N}$, so only remainder estimates have to be checked.
The condition $\sigma +\delta \notin {\Bbb Z}$ is imposed in order to
make \cite{G4, Lemma 3.2} applicable without an $\varepsilon $-reservation.

The coefficients found for
different rays in $V$ coincide with those for the ray $\Bbb R_-$ in view of
the holomorphy (as in \cite{GS1, Lemma 2.3}).
\qed\enddemo

\proclaim{Theorem 2.2} Let $\Cal S^{(N)}_\lambda $ be defined as
in Theorem {\rm 2.1} (with the stated regularity properties and estimates) 
for each each $N=1,2,\dots$, and
assume moreover that $\Cal S^{(N)}_\lambda =\deln \Cal 
S^{(1)}_\lambda $ for all $N$.
Then for $N>(\sigma +n-1)/m$, 
$$
\tilde {\frak s}^{(N)}_{\sigma +n-1 }(x')=
\int_{\Bbb R^{n-1}}{\frak s}^{(N)h}_{ -mN-n+1}(x',\xi
',-1)\,\d\xi'=
\int_{\Bbb R^{n-1}}{\frak s}^{(1)h}_{ -m-n+1}(x',\xi ',-1)\,\d\xi'.\tag2.10
$$

Define the ``log-transform'' $S$ of $\Cal S^{(1)}_\lambda $ as an
operator whose symbol $s(x',\xi ')$ for $|\xi '|\ge 1$ is deduced from
the symbol
$\frak s^{(1)}(x',\xi 
',\lambda )$ of $ {\Cal S}^{(1)}_\lambda $ by $$
s(x',\xi ')=\tfrac i{2\pi }\int_{\Cal C}\log\lambda \, \frak
s^{(1)}(x',\xi ',\lambda )\, d\lambda ,\quad s_{\sigma -j}(x',\xi ')=\tfrac i{2\pi }\int_{\Cal C}\log\lambda \, \frak
s^{(1)}_{\sigma -m-j}(x',\xi ',\lambda )\, d\lambda 
\tag2.11
$$
(with a curve $\Cal C$ in ${\Bbb C}\setminus \crm$ encircling ${\Bbb
C}\setminus V'$). It is 
a classical $\psi $do symbol of order $\sigma $, and for $N>(\sigma +n-1)/m$,  
$$
\aligned
\tilde{\frak s}^{(N)}_{\sigma +n-1
}(x')
&=-\tfrac1m\int_{|\xi '|=1}s_{1-n}(x',\xi ')\,\d S(\xi '),\text{ hence}\\
\tr\tilde{\frak s}^{(N)}_{\sigma +n-1
}(x')&=-\tfrac1m\operatorname{res}_{x'}S.
\endaligned\tag 2.12
$$

\endproclaim

\demo{Proof}
The first equality in (2.10) repeats (2.5), and the last equality follows
in the way explained in \cite{G4, Rem.\ 3.12} from the fact that
$\frak s^{(N)}(x',\xi ',\lambda )
=\tfrac{\partial  _\lambda
^{N-1}}{(N-1)!}\frak s^{(1)}(x',\xi ',\lambda ) $. 

Now consider the ``log-transform''; the integrability in $\lambda $ is
assured by the second line in (2.3). The verification that $s$ is a
classical $\psi $do symbol of order $\sigma $ with homogeneous terms
$s_{\sigma -j}$ goes exactly as in
\cite{G4, Th.  4.5, (4.40)--(4.42)ff.} (with $\sigma +\sigma '$
replaced by $\sigma $, $\frac14$ replaced by $\delta $);
we shall spare a repetition of details. The first formula in (2.12)
represents the fact that
$$
\int_{\Bbb R^{n-1}}{\frak s}^{(1)h}_{ -m-n+1}(x',\xi ',-1)\,\d\xi'
=-\tfrac1m\int_{|\xi '|=1}\tfrac i{2\pi }\int_{\Cal C}\log\lambda
\,{\frak s}^{(1)h}_{ -m-n+1}(x',\xi ',\lambda )\,d\lambda \d S(\xi'),
$$
where the log-integral is turned into an integral along ${\Bbb R}_-$
and the homogeneity is used in the application of polar coordinates, 
by \cite{G4,
Lemmas 1.2 and 1.3} in dimension $n-1$. 
Taking the fiber trace, we get the second formula in (2.12) by
definition, using (1.19) in dimension $n-1$.
\qed 
\enddemo

We shall now see how these theorems apply to show that 
$BG^{(N)}_\lambda $ has a trace expansion
(1.23) with $a'_0=0$ and $a''_0$ equal to a residue. In
\cite{GSc1} the
special case where $G_\lambda $ is
defined from a principally scalar Laplacian with Dirichlet condition
was studied (giving a full
expansion in log-powers of orders going to $-\infty $);  we now
allow more general $\{P_1,T\}$ and just show the expansion down to and
including the crucial term with $(-\lambda )^{-N}$. 
 
The order of $P_1$ is here a positive integer $m$. $B$ is of the form $B=P_++G$, of order $\sigma $ and with $G$ of 
class 0; here $\sigma \in{\Bbb Z}$
when $P\ne 0$, and $\sigma $ can be any real number when $P=0$.

In order to apply Theorems 2.1 and 2.2 we shall show that $G_\lambda $ can be rewritten
in a form that shows a better fall-off in $\lambda $, at the cost of
augmenting the order of the operator. 

\proclaim{Lemma 2.3} 

$1^\circ$ 
The s.g.o.\ $G_\lambda $ and its derivatives $G^{(N)}_\lambda =\deln
G_\lambda $ may be written in the form$$\aligned
G_\lambda &=\lambda ^{-1}P_1G_\lambda,\\
G^{(N)}_\lambda &=(-\lambda )^{-N}P_1 G^{(N)\prime}_\lambda ,
\endaligned\tag 2.13$$
 where $G^{(N)\prime}_\lambda $, described in the proof below, is a
singular Green operator of order $-m$, class $0$ and regularity $+\infty
$.

$2^\circ$ For each $N\ge 1$, the s.g.o.\ $BG^{(N)}_{\lambda }=(P_++G)G^{(N)}_{\lambda }$ has
order $\sigma -mN$, class $0$ and regularity $\sigma $. Moreover, it can be
written as $$
BG^{(N)}_{\lambda }=(-\lambda )^{-N}BP_1G^{(N)\prime}_\lambda ,
\tag2.14
$$
where $BP_1G^{(N)\prime}_\lambda  $ has order $\sigma $, class
$0$ and regularity $\sigma +\frac12$. 

$3^\circ$ In local coordinates, the normal trace $$
 {\Cal U}^{(N)}_\lambda =\tr_n(BG^{(N)}_\lambda)=\tfrac{\partial  _\lambda
^{N-1}}{(N-1)!}\Cal
U^{(1)}_\lambda 
\tag2.15
$$
is a $\psi $do family on $\Bbb R^{n-1}$
 of order $\sigma -mN$ and regularity
$\sigma -\frac14$. It can also be
written as $$
 {\Cal U}^{(N)}_\lambda =(-\lambda )^{-N} {\Cal U}^{(N)\prime}_\lambda
,\quad  {\Cal U}^{(N)\prime}_\lambda =
\tr_n(BP_1G^{(N)\prime}_\lambda ),\tag2.16
$$
where $ {\Cal U}^{(N)\prime}_\lambda $ has order $\sigma $ and regularity
$\sigma +\frac14$. 

$4^\circ$ 
The symbol $\frak u^{(N)}(x',\xi ',\lambda )$ of $\Cal
U^{(N)}_\lambda $, with the expansion in (quasi-)homogeneous terms $\frak
u^{(N)}(x',\xi 
',\lambda )\sim \sum_{j\in \Bbb N}\frak u^{(N)}_{\sigma -mN-j}(x',\xi
',\lambda )$ (the $\frak u^{(N)}_{r}$ being homogeneous of degree $r$ in $(\xi ',\mu )$
on each ray 
$ \lambda =-\mu ^me^{i\theta
}$, $\mu >0$), satisfies:
$$
\aligned
|\partial _{x',\xi '}^{\beta ,\alpha } [\frak u^{(N)}(x',\xi ',\lambda
)-\sum_{j<J}\frak u^{(N)}_{\sigma   -mN-j}(x',\xi ',\lambda )]|
&\leg \ang{\xi '}^{\sigma   - \frac14  -|\alpha |-J}\ang{\xi ',\mu
}^{-mN+\frac14 },\\
|\partial _{x',\xi '}^{\beta ,\alpha } [\frak u^{(N)}(x',\xi ',\lambda
)-\sum_{j<J}\frak u^{(N)}_{\sigma   -mN-j}(x',\xi ',\lambda )]|
&\leg \ang{\xi '}^{\sigma   +\frac14 -|\alpha |-J}\ang{\xi ',\mu
}^{-\frac14 }\mu ^{-mN},
\endaligned
\tag2.17
$$
on the rays in $V$, for all $\alpha 
,\beta ,J$. 
\endproclaim

\demo{Proof} We start by noting that, since $P_1$ is a differential operator,
$$
\aligned
Q_\lambda +\lambda ^{-1}&=Q_\lambda +\lambda ^{-1}(P_1-\lambda
)Q_\lambda =\lambda ^{-1}P_1Q_\lambda \text{ on
}\widetilde X,\\ 
R_\lambda +\lambda ^{-1}&=R_\lambda +\lambda
^{-1}(P_1-\lambda )R_\lambda =\lambda ^{-1}P_1(Q_{\lambda ,+}+G_\lambda )\\
&=\lambda ^{-1}[(P_1Q_\lambda )_++P_1G_\lambda ]=Q_{\lambda
,+}+\lambda ^{-1}
+\lambda ^{-1}P_1G_\lambda \text{ on }X,
\endaligned\tag2.18
$$ 
which implies the first formula in (2.13). For the second formula, we
calculate:
$$\aligned
G^{(N)}_\lambda &=P_1\deln(\lambda ^{-1}G_\lambda )
=P_1\sum_{j=0}^{N-1}c_j(-\lambda )^{-1-j}G^{(N-j)}_\lambda \\
&=(-\lambda )^{-N}P_1\sum_{j=0}^{N-1}c_j
((P_1-\lambda )-P_1)^{N-1-j}G^{(N-j)}_\lambda \\
&=(-\lambda )^{-N}P_1\sum_{j=0}^{N-1}c_j\sum  _{k=0}^{N-1-j}c'_kP_1^k(P_1-\lambda )^{N-1-j-k}G^{(N-j)}_\lambda ;
\endaligned\tag2.19$$ 
here the sums over $j$ and $k$ define an s.g.o.\
$G^{(N)\prime}_\lambda $ of order $-m$, class $0$ and
regularity $+\infty $.

For $2^\circ$ we use the rules in \cite{G1, Prop.\ 2.3.14}, which give
that $P_+$ and $G$ have order and regularity $\sigma $, whereas
$P_+P_{1,+}=(PP_1)_+-L(P,P_1)$ has order $\sigma +m$ with $(PP_1)_+$ of regularity $\sigma +m$,
$L(P,P_1)$ an s.g.o.\ of class $m$ and hence regularity $\sigma
+\frac12$, and likewise $GP_{1,+}$ of order $\sigma +m$ and class $m$ and hence regularity $\sigma
+\frac12$. The compositions with $G^{(N)}_\lambda $ resp.\
$G^{(N)\prime}_\lambda $ inherit these regularities by the rules in
\cite{G1, Sect.\ 2.7}. 

For $3^\circ$,  there is a loss of $\frac14$ in the
regularity when the general rule of \cite{G4, Lemma 3.4} is
applied.  

The information in $4^\circ$ follows from the definition of the class
of symbols of the stated regularity,  when $\sigma <0$. When $\sigma
\ge 0$, the regularity information itself gives weaker estimates when
$\sigma -|\alpha |-J\ge 0$. But here we can use the device introduced
in the proof of \cite{G4, Prop.\ 4.3}: Compose $\Cal U^{(N)}_\lambda $ to the
left with $\Lambda ^\varrho \Lambda ^{-\varrho }$, $\Lambda
=\operatorname{OP}(\ang{\xi '})$, with $\varrho \ge \sigma +1$. Taking
$\Lambda ^{-\varrho }$ together with $B$ one finds that
$\Lambda ^{-\varrho }\Cal U^{(N)}_\lambda $ satisfies the regularity
statements with $\sigma $ replaced by $\sigma -\varrho $, hence the
estimates with the same replacement, and the desired estimates follow
after composition with $\Lambda ^\varrho $.
\qed\enddemo

Theorems 2.1 and 2.2 apply straightforwardly to the symbols described
in this lemma:

\proclaim{Theorem 2.4} Let ${\Cal U}^{(N)}_\lambda $ be as in Lemma
{\rm 2.3}. When $N>(\sigma +n-1)/m$, the kernel of 
$ {\Cal
 U}^{(N)}_\lambda $ on the diagonal has an expansion$$
K( {\Cal U}^{(N)}_\lambda ,x',x')=\sum_{0\le l\le \sigma +n-1}\tilde
{\frak u}^{(N)}_l(x')(-\lambda )
^{\frac{\sigma +n-1 -l}m -N}+O(\lambda ^{-N-\frac1{4 m}(+\varepsilon )}),
\tag2.20
$$
where $(+\varepsilon )$ indicates the addition of an $\varepsilon >0$
if $\sigma +\frac14\in {\Bbb Z}$ (no addition otherwise).
The coefficient of $(-\lambda ) ^{-N}$ is$$
\tilde {\frak u}^{(N)}_{\sigma +n-1}(x')=
\int_{\Bbb R^{n-1}}{\frak u}^{(N)h}_{ -mN-n+1}(x',\xi
',-1)\,\d\xi'=
\int_{\Bbb R^{n-1}}{\frak u}^{(1)h}_{ -m-n+1}(x',\xi ',-1)\,\d\xi'\tag2.21
$$
if $\sigma +n-1\in{\Bbb N}$, zero if $\sigma +n-1\notin{\Bbb N}$.

 The ``log-transform'' $U$, defined as an operator with symbol
$u(x',\xi ')$ deduced for $|\xi '|\ge 1$ from the symbol
$\frak u^{(1)}(x',\xi 
',\lambda )$ of $ {\Cal U}^{(1)}_\lambda $ by $$
u(x',\xi ')=\tfrac i{2\pi }\int_{\Cal C}\log\lambda \, \frak
u^{(1)}(x',\xi ',\lambda )\, d\lambda, \quad
u_{\sigma -j}(x',\xi ')=\tfrac i{2\pi }\int_{\Cal C}\log\lambda \, \frak
u^{(1)}_{\sigma -m-j}(x',\xi ',\lambda )\, d\lambda 
\tag2.22
$$
(and extended smoothly for $|\xi '|\le 1$), is a classical $\psi $do
of order $\sigma $ such that the coefficient of $(-\lambda )^{-N}$ 
in {\rm (2.20)} satisfies:$$
\aligned
\tilde{\frak u}^{(N)}_{\sigma +n-1
}(x')&=-\tfrac1m\int_{|\xi '|=1}u_{1-n}(x',\xi ')\,\d S(\xi '),\\
\tr\tilde{\frak u}^{(N)}_{\sigma +n-1
}(x')&=-\tfrac1m\operatorname{res}_{x'}U;
\endaligned
\tag 2.23
$$
here $u_{1-n}$ is defined as zero if $\sigma +n-1\notin{\Bbb N}$.
\endproclaim 

\demo{Proof} When $\sigma \in{\Bbb Z}$ (which must hold if $P\ne 0$),
we apply Theorems 2.1 and 2.2 with $\delta =\frac14$. When $\sigma \notin{\Bbb
Z}$,
we apply the theorems with $\delta =\frac14$ if $\sigma
+\frac14\notin{\Bbb Z}$,  $\delta <\frac14$ if $\sigma
+\frac14\in{\Bbb Z}$.\qed  
\enddemo

Integrating over the local coordinate patches and carrying the pieces
back to the manifold, we find:

\proclaim{Corollary 2.5} When $N>(\sigma +n-1)/m $, $\Tr(BG^{(N)}_\lambda )$
has an expansion
$$
\Tr(BG^{(N)}_\lambda )=\sum_{0\le l\le \sigma +n-1}\tilde
{\frak u}^{(N)}_l(-\lambda )
^{\frac{\sigma +n-1 -l}m -N}+O(\lambda ^{-N-\frac1{4 m}(+\varepsilon )})
\tag2.24
$$
(with $(+\varepsilon )$ understood as in Theorem {\rm 2.4}),
where the coefficient of $(-\lambda ) ^{-N}$ is$$
\tilde {\frak u}^{(N)}_{\sigma +n-1}=
\int_{X'}\tr\tilde{\frak u}^{(N)}_{\sigma +n-1}(x')\,dx'=
-\tfrac1m\int_{X'}\operatorname{res}_{x'} U\,dx'=-\tfrac1m
\operatorname{res} U,
$$
if $\sigma +n-1\in{\Bbb N}$, zero otherwise.
The expansion may also be written as:
$$
\Tr(BG^{(N)}_\lambda )=\sum_{1\le j\le \sigma +n}b^{(N)}_j(-\lambda )
^{\frac{\sigma +n-j}m -N}+O(\lambda ^{-N-\frac1{4 m}(+\varepsilon )}),
\tag2.25
$$
with $b^{(N)}_{\sigma +n}=-\frac1m \operatorname{res}U$; in other words,
$$
l_0(BG_\lambda )=-\tfrac1m \res U.\tag 2.26
$$

\endproclaim

\demo{Proof} We use here that $G_\lambda ^{(N)}$ cut down to interior
coordinate patches is of order $-\infty $ and rapidly decreasing in
$\lambda $, so that such coordinate patches contribute only to the
remainder in (2.24). The alternative formulation in (2.25)
is obtained by denoting $l+1=j$ and relabelling
$$
\tilde
{\frak u}^{(N)}_{l}=b ^{(N)}_{l+1}. \quad\square
$$
\enddemo

In the case where $m>\sigma +n$, this result can simply be added to
(1.4), implying the validity of (1.3) with $N=1$ in this case.

We also want to establish (1.3) for general $N$, $m=2$, which requires
proving a version of \cite{G4, Th.\ 3.6} where $B(Q_{1,\lambda
}-Q_{2,\lambda })_+$ for auxiliary operators of order $m>\sigma +n$ is
replaced by $B(Q^N_{1,\lambda }-Q^N_{2,\lambda })_+$ for auxiliary
operators of order 2. This can be done with the same methods as used for   
\cite{G4, Th.\ 3.6, Th.\ 3.10}, and could also be done more mechanically with 
the technology of Section 4
there. Since there are no new difficulties in this, the explanation
will be brief. Without extra effort, we can let the auxiliary
differential operators have an arbitrary positive order $m$.

\proclaim{Theorem 2.6} Let $P_1$ and $P_2$ be auxiliary elliptic
differential operators on $\widetilde X$ of order $m>0$ with no eigenvalues on ${\Bbb
R}_-$, let $Q_{i,\lambda }=(P_i-\lambda )^{-1}$, and consider
$B(Q^N_{1,\lambda }-Q^N_{2,\lambda })_+$ on $X$,
decomposed in its $\psi $do part and s.g.o.\ part
$$
B(Q^N_{1,\lambda }-Q^N_{2,\lambda })_+=(P\Cal Q^{(N)}_\lambda )_++\Cal
G^{(N)}_\lambda ,
\text{ where }\Cal Q^{(N)}_\lambda =Q^N_{1,\lambda }-Q^N_{2,\lambda }.\tag2.27
$$
For $N>(\sigma +n)/m$, the $\psi $do part and s.g.o.\ parts have trace expansions
$$
\aligned
\Tr(P\Cal Q^{(N)}_{\lambda })_+
&=\sum_{0\le j \le n+\sigma }c^{(N)}_{j}\,(-\lambda )
^{\frac{n+\sigma -j}m-N}+O(\lambda  ^{-N-\frac1{4m}(+\varepsilon ) }),\\
\Tr\Cal G^{(N)}_{\lambda }
&=\sum_{1\le j \le n+\sigma }d^{(N)}_{j}\,(-\lambda )
^{\frac{n+\sigma -j}m-N}+O(\lambda  ^{-N-\frac 1{4m}(+\varepsilon ) }),
\endaligned
\tag2.28
$$ 
with $(+\varepsilon )$ understood as in Theorem {\rm 2.4}.
Here the coefficients of $(-\lambda )^{-N}$ are zero if
$\sigma +n\notin{\Bbb N}$, and otherwise, in terms of local coordinates, 
$$
c^{(N)}_{\sigma +n}=-\tfrac1m \int_{X}\operatorname{res}_x(P(\log
P_1-\log P_2))\,dx,\tag2.29
$$
and 
$$
d^{(N)}_{\sigma
+n}=-\tfrac1m\int_{X'}\operatorname{res}_{x'}S\,dx',\text{ with }
S=\tfrac i{2\pi }\int_{\Cal C}\log\lambda \, \tr_n(\Cal G^{(1)}_\lambda
)\,d\lambda .\tag 2.30
$$
 
If $m$ is even, {\rm (1.12)} holds with the residue defined in {\rm \cite{FGLS}}.
\endproclaim

\demo{Proof} The $\psi $do part is dealt with by methods as in \cite{G4,
Sect.\ 2}. The crucial fact is that the symbol of $\Cal Q^{(N)}_\lambda
$, as a difference between two iterated resolvent symbols, has
homogeneous terms that are rational functions of $(\xi ,\lambda )$
with a one step better fall-off in $\lambda $ than the individual
symbols of the iterated  resolvents $Q^N_{i,\lambda
}$, as in \cite{G4, Prop.\ 2.1} (it is only the leading term that
needs some thought). For the composition with $P$ this
implies that there is a strictly homogeneous symbol, integrable at $\xi
=0$, which produces the term with 
$(-\lambda )^{-N}$ in the first line of (2.28). Then a reformulation 
in polar coordinates, plus a comparison
with log-formulas as in the proof of \cite{G4, Th. 2.2}, lead to a
diagonal kernel expansion that integrates over $X$ to give (2.28)
with (2.29).

For the s.g.o.\ term, if $\sigma \in{\Bbb Z}$, the symbolic properties of $\Cal Q^{(N)}_\lambda
$ are used as in \cite{G4, Th.\ 3.6} to see that there is a strictly
homogeneous term in the symbol of $\tr_n\Cal G^{(N)}_\lambda $, 
integrable at $\xi '
=0$,  that produces the term with 
$(-\lambda )^{-N}$ in the second line of (2.28).
It is interpreted as in \cite{G4, Th. 3.10, Rem.\ 3.12} to give
(2.30). If $\sigma \notin {\Bbb Z}$, there is no such term, only the
remainder. 

When $m$ is even, $\log P_1-\log P_2$ has the transmission property,
and 
$$
G'=-L(P,\log P_1-\log P_2)
+G(\log P_1-\log P_2)_+
$$
is a singular Green operator in the calculus. It is seen as in the
proof of \cite{G4, Th.\ 3.6} that $\res_{x'}S=\res_{x'}\tr_n G'$ (in
the relevant term, the log-integration can be moved outside
$\tr_n$). Then indeed,
$$
\aligned
C_0(B,P_{1,+})-C_0(B,P_{2,+})
&=-\tfrac 1m \operatorname{res}(P(\log
P_1-\log P_2))_+
-\tfrac 1m \operatorname{res}_{X'}S\\
&=-\tfrac 1m \operatorname{res}(B(\log
P_1-\log P_2)_+).\quad\square
\endaligned\tag 2.31
$$
\enddemo

As remarked earlier, the interpretation in \cite{G4, Th.\ 3.6} of the
trace defect as a 
residue in the sense of \cite{FGLS} holds only when $m$ is even;
otherwise it can be regarded as a residue defined in a more general
sense. We note in passing that there is a misprint in the statement
of the theorem there; $\sigma -\frac14$ should be replaced by
$\sigma +\frac14$ in line 6 of page 1691.

\proclaim{Corollary 2.7} $1^\circ$ The expansion {\rm (1.3)} holds
with $N=1$ for general
systems $\{P_1,T\}$ of order $m>\sigma +n$. 

$2^\circ$  For general second-order elliptic operators $P_1$, the
operator family $BQ^N_{\lambda ,+}$ has trace expansions {\rm (1.23)}
 when $N>(\sigma +n)/2$.

$3^\circ$ The expansion {\rm (1.3)} holds for general
systems $\{P_1,T\}$ of order $2$, with $N>(\sigma +n)/2$. 
\endproclaim

\demo{Proof} For (1.2) we have already observed that it follows by
adding the result of Corollary 2.5 with $m>\sigma +n$, $N=1$, to (1.4).

For $2^\circ$, we consider a general operator $P_1$ together with a
special choice $P_2$ as in \cite{GSc1}; then the result from there on
the expansion of  $\Tr(B(P_2-\lambda )^{-N}_+)$ taken together with Theorem
2.6 on $\Tr(B((P_1-\lambda )_+^{-N}-(P_2-\lambda )_+^{-N}))$ implies the statement.

Now (1.3) is obtained for a general second-order realization 
$P_{1,T}$ by combination with
Corollary 2.5 for $m=2$, $N>(\sigma +n)/2$.\qed
\enddemo

There is more information on the case of general even $m$
in Section 6. 

\example{Remark 2.8} For auxiliary operators $P_1$ of order 1,
Theorems 2.4--2.6 are valid. Results as in Corollary 2.7 can also be 
shown, under slightly restrictive circumstances:

Assume that $P_1$ can be chosen skew-selfadjoint of order 1, acting in the
bundle $\widetilde E$ over $\widetilde
X$. Then we can use the trick (found e.g.\ in \cite{GS1}) of
introducing a ``doubled''
operator (also skew-selfadjoint)
$$
\Cal P_1=\pmatrix 0&P_1\\P_1&0\endpmatrix,\tag2.32
$$
whose resolvent is 
$$
(\Cal P_1-\lambda )^{-1}=\pmatrix \lambda (P_1^2-\lambda ^2)^{-1}&
P_1 (P_1^2-\lambda ^2)^{-1}\\ P_1 (P_1^2-\lambda ^2)^{-1}&\lambda (P_1^2-\lambda ^2)^{-1}
\endpmatrix,\tag2.33
$$
for $\lambda \in {\Bbb C}\setminus i{\Bbb R}$.
Define also 
$$
\Cal B=\pmatrix B&0\\0&B\endpmatrix,\tag2.34
$$
then if $B$ is of so low order that the trace exists ($\sigma +n<1$),
$$
\Tr(\Cal B(\Cal P_1-\lambda )^{-1}_+)=2\lambda \Tr (B(P_1^2-\lambda ^2)^{-1}_+).
$$
Since $P_1^2$ is selfadjoint negative, the right-hand side has an
expansion in powers of $\lambda $ on suitable rays, by the
preceding results. For the
left-hand side, this gives a trace expansion for one special case of
a first-order operator. Thanks to Theorem 2.6, we can then also get
expansions for $\Tr(\Cal B (\Cal P'-\lambda )^{-1}_+)$ for other
choices of first-order operators $\Cal P'$. Take e.g.
$$
\Cal P'=\pmatrix P_1&0\\0&P_1\endpmatrix,\tag2.35
$$
then 
$$
\Tr(\Cal B(\Cal P'-\lambda )^{-1}_+)=2 \Tr (B(P_1-\lambda )^{-1}_+),
$$
from which we obtain an expansion of $\Tr (B(P_1-\lambda
)^{-1}_+)$ as in (1.2),  on suitable rays. 
Now Theorem 2.6 can be
used again, to allow other first-order elliptic operators $P_2$ in the
place of $P_1$.

If $B$ is not of low order, one must work with derivatives $\deln$ of
$\Tr(\Cal B(\Cal P_1-\lambda )^{-1}_+)$ to get the result, which gives
more complicated calculations.

Thanks to Corollary 2.5, the analysis extends to cases where
$(P_2-\lambda )_+^{-1}$ is replaced by $(P_{2,T}-\lambda )^{-1}$,
where $P_2$ is provided 
with an elliptic differential boundary condition $Tu=0$. 

Let us remark
that not all first-order elliptic operators have elliptic differential
boundary conditions. For Dirac operators there are in any case the
{\it pseudodifferential} Atiyah-Patodi-Singer boundary conditions (and their
generalizations), but the present set-up would not include
them; this requires different efforts (cf.\ e.g.\ [G2]), since the 
singular Green part
$G_\lambda $ then has low regularity.
(There is much other current literature on that subject, also from the
author's side, but we shall not lengthen the reference list with this.) 
\endexample

Having established (1.3) for general $N$, $m=2$, we go on to the
detailed analysis of $C_0(B, P_{1,T})$, treating the
contributions from the five terms in (1.22) one by one.

\subhead 3. The constant coming from $P_+Q_{\lambda ,+}$; localization
\endsubhead 

We now begin the analysis of the contributions to $a''_0=C_0(B, P_{1,T})$.
For the contributions from $BQ^{N}_{\lambda ,+}=P_+Q^{N}_{\lambda ,+}+GQ^{N}_{\lambda ,+}$ 
in (1.22) (which together give the term $C_0(B,P_{1,+})$), it is an advantage to start with a very simple case where 
$P_1$ is similar to the Laplacian.

As in \cite{G5}, we make a reduction to local coordinates and choose
$P_1$ with essentially constant coefficients there.
This goes as follows: $X$ is covered by a system of open subsets $U_j$,
$j=1,\dots,J$, of $\widetilde X$, with trivializations $\Phi
_j\:\widetilde E|_{U_j}\to V_j\times \Bbb C^M$ ($M=\dim \widetilde
E$) and base maps $\kappa _j\:U_j\to V_j$, with $V_j$ bounded in
$\Bbb R^n$. Now $\{\psi _j\}_{1\le
j\le J}$ is an 
associated partition of unity (with $\psi _j\in C_0^\infty (U_j)$),
and $\varphi _j\in C_0^\infty (U_j)$ with $\varphi _j=1$ on
$\operatorname{supp}\psi _j$. Then 
$$
B=\sum_{1\le j\le J}\psi _jB= \sum_{1\le j\le J}\psi _jB\varphi
_j+\sum_{1\le j\le J}\psi _jB(1-\varphi _j),
$$
where the last sum is of order $-\infty $; for this part
it is well-known that \linebreak$C_0(\psi _jB(1-\varphi _j),P_{1,+})=\Tr (\psi _jB(1-\varphi _j))$.
So it remains to treat each of the terms $\psi _jB\varphi
_j$. 

For the present case of a manifold with boundary we can assume that,
say, the sets with $j=1,\dots,J_0$ intersect $\partial X$ and the
remaining sets with $j=J_0+1,\dots, J$ have closures lying in the interior of 
$X$, and we can for $j\le J_0$ take the $V_j$ of the form 
$W_j\times \,]-c,c[\,$, $W_j\subset \Bbb R^{n-1}$ such that $W_j=\kappa
_j(U_j\cap \partial X)$.

For $j>J_0$, $\psi _jB\varphi _j=\psi _j(P+\Cal R_j)\varphi _j$ with
$\Cal R_j$ a $\psi $do of order $=\infty $ (since s.g.o.s are
smoothing on the interior of $X$), so these terms are essentially
covered by the analysis in \cite{G5}, and by the account we give for
the $\psi $do part below. Therefore we can restrict the attention to
one of the terms with $j\le J_0$, say, $\psi _1B\varphi _1$.  

As noted in \cite{G5},  one can assume that $X$ is already covered by
$U_{10}, U_2,\dots ,U_J$ with $\overline U_{10}$ compact in $U_1$,
$\psi _1$ and $\varphi _1$ supported in $U_{10}$, and introduce
$U_1'=U_1$, $U'_j=U_j\setminus \overline U_{10}$ for $j\ge 2$, as a
new cover of $X$ with associated partition of unity $\psi '_j$ and
functions $\varphi '_j\in C_0^\infty (U'_j)$ equal to 1 on $\supp
\psi '_j$, $1\le j\le J$. This allows us to choose 
$$
P_1u=\sum_{1\le j\le J}\varphi '_j[(-\Delta )I_M((\psi '_ju)\circ {\Phi
^*_j}^{-1})]\circ \Phi _j^*,\tag3.1
$$ 
as the auxiliary elliptic operator, such that when 
$\psi _jB\varphi
_j$ is carried over to $\widetilde B$ in the coordinate system $V_1$,
$\psi _jB\varphi
_j(P_1-\lambda )^{-1}$ is carried over to $\widetilde B(-\Delta
-\lambda )^{-1}I_M+\Cal R_\lambda $ in $V_1$, where $\Cal R_\lambda $
is of order $-\infty $ with a trace that is $O(\lambda ^{-N'})$, any
$N'$
(with similar properties of $\lambda $-derivatives). $I_M$ stands for
the $M\times M$ identity matrix, not mentioned explicitly in the following.
(Cf.\ \cite{G5}
for more details.)

This reduces the problem to a calculation of the trace expansion of
operators of the form $\widetilde B(-\Delta -\lambda )^{-N}$ on $\rnp$,
where we can work out the kernel explicitly. (In \cite{G5}, $1-\Delta
$ was used rather than $-\Delta $, since $\log(1-\Delta )$ makes sense
on $\rn$, but the lack of homogeneity is a disadvantage when we
calculate more complicated boundary contributions; instead we shall modify
the symbol near 0 when needed.)

After the reduction, we again write $\widetilde B$ as $P_++G$, with symbols
$p(x,\xi )$ 
resp.\ $g(x',\xi ',\xi _n,\eta _n)$, functions on $\Bbb R^{2n}$.

In the rest of the present section we consider the contributions from
$P_+Q_{\lambda ,+}$ and its iterated versions; recall that the order $\sigma $ of $P$ is an integer.
The operator $P_+Q^N_{\lambda ,+}$ is written as a sum of two terms$$
P_+Q^N_{\lambda ,+}=(PQ^N_{\lambda })_+-L(P,Q^N_{\lambda
})
=\tfrac{\partial  _\lambda ^{N-1}}{(N-1)!}[(PQ_{\lambda
})_+-L(P,Q_{\lambda })],
\tag3.2
$$
that are treated in different ways. The first term is the truncation to
$\Bbb R^n_+$ of the $\psi $do $PQ_\lambda ^N$ on $\Bbb R^n$, whereas
the second term, the ``leftover term'',  is a composition of s.g.o.s;
see (3.16)ff.\ later. We see in the following that for the first term one
can use the method of \cite{G5}, for the second term that of Section 2.

For the first term, the pointwise calculations are essentially as in
\cite{G5}, relative to the closed manifold $\widetilde X$, and we get
the contribution to the basic zeta value as the
integral over $X$ of the pointwise contribution (pulled back from
local coordinates). A short explanation will now be given.

\proclaim{Proposition 3.1} When $N>(\sigma +n)/2$, the kernel of
$PQ^N_{\lambda }$ on the 
diagonal $\{x=y\}$ has an expansion:$$
\aligned
K(PQ^N_{\lambda },x,x)&=\sum_{0\le j<\sigma +n}(-\lambda
)^{\frac{\sigma +n-j}2-N}c_j(x)\\
&\quad +\tfrac{\partial  _\lambda
^{N-1}}{(N-1)!}[(-\lambda )^{-1}\log(-\lambda )c'_0(x)+(-\lambda
)^{-1}c''_0(x)]+ O(\lambda ^{-N-\frac12+\varepsilon })\\
&=\sum_{0\le j<\sigma +n}(-\lambda
)^{\frac{\sigma +n-j}2-N}c_j(x)+(-\lambda )^{-N}\log(-\lambda
)c'_0(x)\\
&\quad+(-\lambda
)^{-N}(c''_0(x)-\alpha _Nc'_0(x))+ O(\lambda ^{-N-\frac12+\varepsilon }),
\endaligned
\tag3.3
$$ 
any $\varepsilon >0$, for $\lambda \to\infty $ on rays in $\Bbb
C\setminus \Bbb R_+$. Here $\alpha
_N=\sum_{1\le k <N}\frac1k$, and$$
\aligned
c_j(x)&=\int_{\Bbb R^n}p^h_{\sigma -j}(x,\xi )(|\xi
|^2+1)^{-N}\,\d \xi,\\
c'_0(x)&=\tfrac12\int_{|\xi |=1}p_{-n}(x,\xi)\,\d  S(\xi),\\
c''_0(x)&=\slint p(x,\xi )\,\d \xi;  
\endaligned\tag3.4$$
then (cf.\ {\rm (1.16), (1.19)})
$$
\tr c'_0(x)=\tfrac12\operatorname{res}_xP,\quad \tr c''_0(x)= \TR_x P.\tag3.5
$$

\endproclaim 
\demo{Proof} This is a purely pseudodifferential situation, and the proof goes as in \cite{G5}, except that we must in
addition account for the
effect of the $\lambda $-derivative. Using the expansion
$p\sim\sum_{j\in\Bbb N}p_{\sigma -j}(x,\xi )$ in terms $p_{\sigma
-j}(x,\xi )$ homogeneous of degree $\sigma -j$ in $\xi $ 
for $|\xi |\ge 1$, we write$$
p(x,\xi )=\sum_{0\le j< \sigma +n}p_{\sigma -j}(x,\xi )
+p_{-n}(x,\xi )+p_{<-n}(x,\xi ).
$$
The finite-part integral $\tslint p(x,\xi  )\,\d\xi $ is defined as
recalled above after (1.16).
Here (cf.\ e.g.\ \cite{G3, (1.18)})
$$
\aligned
\slint p_{\sigma -j}(x,\xi )\,\d \xi &=
\int_{|\xi |\le 1}(p_{\sigma -j}(x,\xi
)-p^h_{\sigma -j}(x,\xi ))\,\d \xi, \text{ for }\sigma -j>-n,\\
\slint p_{-n}(x,\xi )\,\d \xi &=\int_{|\xi |\le 1}p_{-n}(x,\xi
)\,\d \xi,\\
\slint p_{<-n}(x,\xi )\,\d  \xi&=\int_{\Bbb R^n}p_{<-n}(x,\xi
)\,\d  \xi 
,\endaligned
\tag3.6
$$
where  $p^h_{\sigma -j}$ denotes the strictly
homogeneous version of $p_{\sigma -j}$ (the extension of $p_{\sigma
-j}|_{|\xi |\ge 1}$
into $\{1>|\xi |>0\}$ that is homogeneous for $|\xi |>0$). 

The kernel of $PQ^N_\lambda $ is the integral in $\xi $ of its symbol
times $e^{i(x-y)\cdot \xi }$. Its value for $x=y$ is then 
$$
K(PQ^N_{\lambda },x,x)=K(P(-\Delta -\lambda )^{-N},x,x)=\int_{\Bbb R^n}p(x,\xi )(|\xi
|^2-\lambda )^{-N}\,\d\xi ,\tag3.7
$$
which we shall expand in powers of $\lambda $.

The symbols $p^h_{\sigma -j}$
with $j<\sigma +n$ 
are integrable at $\xi =0$, so these terms $p_{\sigma -j}$ produce
the sum over $j<\sigma +n$ in (3.3) plus $(-\lambda
)^{-N}\tslint\sum_{j<\sigma +n}p_{\sigma -n}\,\d\xi $ by
homogeneity as in the proof of \cite{G5, Lemma 4.1}, with an error
that is $O(\lambda ^{-N-1})$.

For the symbol $p_{-n}$, one has for any $N\ge 1$ that
$$
\int_{\Bbb R^n}p_{-n}(x,\xi )(|\xi |^2-\lambda )^{-N}\,\d \xi 
=\tfrac{\partial  _\lambda
^{N-1}}{(N-1)!}\int_{\Bbb R^n}p_{-n}(x,\xi )(|\xi |^2-\lambda
)^{-1}\,\d \xi ,\tag3.8
$$
and the proof in \cite{G5, Lemma 4.1} can be used
word for word, when it is furthermore remarked that the expansion of
$\log(1+s)$ allows 
taking derivatives; this gives 
$$
\aligned
\int_{\Bbb R^n}&p_{-n}(x,\xi )(|\xi |^2-\lambda )^{-N}\,\d \xi \\
&=\tfrac{\partial  _\lambda
^{N-1}}{(N-1)!}[\tfrac12(-\lambda )^{-1}\log(-\lambda )\int_{|\xi
|=1}p_{-n}(x,\xi )\d S(\xi)
+(-\lambda )^{-1}\slint p_{-n}(x,\xi
)\d\xi ]
+O(\lambda ^{-N-1}) \\
&=(-\lambda )^{-N}\log(-\lambda )c'_0(x)
+(-\lambda )^{-N}[\slint p_{-n}(x,\xi
)\,\d\xi -\alpha _Nc'_0(x)]+O(\lambda ^{-N-1}),
\endaligned
$$
for any $N\ge 1$; $\alpha _N$ comes from differentiating the log-term.

Finally, 
since $p_{<-n}$ and $|\xi |^{1-\varepsilon '}p_{<-n}$ are integrable
for $\varepsilon ' \in \,]0,1]$,
we find $$
\int_{\Bbb R^n}p_{<-n}(x,\xi )(|\xi |^2-\lambda )^{-N}\,\d \xi
=(-\lambda )^{-N}\int_{\Bbb R^n}p_{<-n}(x,\xi )\,\d \xi
+O(\lambda ^{-N-\frac12+\varepsilon }),\tag3.9  
$$
any $N\ge 1$ and $\varepsilon >0$,
by
insertion of an expansion of $(|\xi |^2-\lambda )^{-N}$:
$$
\aligned
(|\xi |^2-\lambda )^{-1}&=
(-\lambda )^{-1}-(-\lambda )^{-1}|\xi |^2(|\xi |^2-\lambda )^{{-1}}
 \\
(|\xi |^2-\lambda )^{-N}&=
(-\lambda )^{-N}-\sum_{0\le k< N}c_k(-\lambda )^{-k}((-\lambda )^{-1}|\xi |^2(|\xi
|^2-\lambda )^{-1})^{N-k}\\
&=
(-\lambda )^{-N}-(-\lambda )^{-N}\sum_{0\le k< N}c_k(|\xi |^2(|\xi
|^2-\lambda )^{-1})^{N-k}\\
&=(-\lambda )^{-N}+O((-\lambda )^{-N-\frac12+\varepsilon }|\xi
|^{1-2\varepsilon }),\; \varepsilon \in \,]0,\tfrac12 [\,.
\endaligned
\tag3.10$$
This shows (3.3). \qed
\enddemo

\proclaim{Proposition 3.2} The operator family $(PQ_\lambda ^N)_+ =(P(-\Delta -\lambda )^{-N})_+$ in {\rm (3.2)} on $\Bbb
R^n_+$ has for $N>(\sigma +n)/2$ expansions as in {\rm (1.23)}, with
$\delta =\frac12-\varepsilon $, any $\varepsilon >0$, for $\lambda \to\infty $ on rays in $\Bbb
C\setminus \Bbb R_+$.
Here $a'_0=\tfrac12\res(P_+)$ and$$
l_0( (P(-\Delta -\lambda )^{-1})_+)=
\int_{\rnp}[\TR_x P-\tfrac12\res_{x,0}(P\log P'_1)]\, dx;
\tag3.11$$
with $\res_{x,0}(P\log P'_1)=0$. 
\endproclaim 

\demo{Proof} This is found by integrating the fiber trace of the
expansion in (3.3) with respect to $x\in\Bbb R^n_+$ (recall that the symbol of $P$ has compact $x$-support in
this situation). This gives rise to the coefficients $a'_0$ and
$l_0$ as described. Here $\res_{x,0}(P\log P'_1)$ is zero, 
because $\log P'_1$ has symbol
$2\log [\xi ]$, so that the log-free term of order 
$-n$ in the symbol of $P\log P'_1$ vanishes for $|\xi |\ge 1$.\qed
\enddemo 

We have included the trivial term with
 $\res_{x,0}(P\log P'_1)$ in (3.11) for the sake of generalizations to other
auxiliary operators $P_2$.

As explained in  the beginning of this section, we can choose $P_1$ on the manifold so that its resolvent is similar
to $(-\Delta -\lambda )^{-1}$ in
the specially selected local coordinates. Then, summing over the coordinate 
patches, we find that $(P(P_1-\lambda )^{-N})_+$ has a trace expansion
(1.23), where $l_0$ is a sum of 
contributions of the form $$
\varphi \,\int [\TR_x
P-\tfrac12\res_{x,0}(P\log P_1)]\,dx\,\psi  ,\tag3.12$$ 
with cutoff functions $\varphi ,\psi $. Also here,
$\res_{x,0}(P\log P_1)$ is 0, since $P_1$
in the local coordinates has the same symbol terms as $P'_1$ for $|\xi
|\ge 1$. Now the formula will be extended to general choices of auxiliary
operator as follows:

When $P_2$ is a general auxiliary elliptic operator of order 2, 
 the kernel of \linebreak$P((P_2-\lambda )^{-N}-(P_1-\lambda )^{-N})$ has an
expansion on the diagonal (calculated in local coordinates):
$$\multline
K(P((P_2-\lambda )^{-N}-(P_1-\lambda )^{-N}),x,x)=\sum_{0\le j \le\sigma +n}s_j^{(N)}(x)(-\lambda
)^{\frac{\sigma +n-j}2-N}\\
+ O(\lambda ^{-N-\frac12 }),\endmultline\tag3.13
$$
with 
$$
s^{(N)}_{\sigma +n}(x)=s_{\sigma +n}(x)=-\tfrac12\int_{|\xi
|=1}\operatorname{symb}_{-n}(P(\log P_2-\log P_1))(x,\xi )\,\d S(\xi ).
$$
The general idea of proof of this is given in \cite{G4, Prop.\ 2.1, Th.\ 2.2},
and the adaptation to the situation with $N$'th powers is given in
\cite{G5, Sect.\ 3}. Since $\log
P_2-\log P_1$ is a classical $\psi $do having the transmission
property (cf.\ \cite{GG, Lemma 2.1}), so is $P(\log P_2-\log P_1)$, so
the integral over $X$ of the fiber trace of (3.14) is precisely
$-\frac12 \res ([P(\log P_2-\log P_1)]_+)$ defined as in
\cite{FGLS}. Thus integration over $X$ of the fiber trace of (3.13)
shows that $[P((P_2-\lambda )^{-N}-(P_1-\lambda )^{-N})]_+$ is as in
Definition 1.2 with $\delta =\frac12$ and
$$\aligned
l_0([P((P_2-\lambda )^{-1}-(P_1-\lambda )^{-1})]_+)
&=-\tfrac12 \res
([P(\log P_2-\log P_1)]_+)\\
&=-\tfrac12 \int_X\res _x(P(\log P_2-\log
P_1))\, dx.
\endaligned\tag3.14
$$
Here the symbol of $P(\log P_2-\log P_1)$ is classical, without
logarithmic terms, so $\res_{x,0}$ of it coincides with $\res _x$ of
it. Then $\res_{x,0}(P\log P_2)=\res_{x}(P(\log P_2 -\log
P_1))+\res_{x,0}(P\log P_1). $ Adding (3.14) to the result for
$l_0((P(P_1-\lambda )^{-1})_+)$, we conclude:

\proclaim{Theorem 3.3}  For a 
$\psi $do $P$ of order $\sigma \in\Bbb Z$ satisfying the transmission
condition 
at $X'$, together with a general elliptic differential operator
$P_2$ of order $2$ having $\rmi$ as a spectral cut, 
the operator family $(P(P_2 -\lambda )^{-N})_+$ has for $N>(\sigma
+n)/2$ expansions as in {\rm (1.23)}, with
$\delta =\frac12-\varepsilon $, any $\varepsilon >0$, for $\lambda
\to\infty $ on rays in the sector $V$ around ${\Bbb R}_-$ where
$p_2^0-\lambda $ is invertible. Here
$$
l_0( (P(P_2 -\lambda )^{-1})_+)=
\int_{X}[\TR_x P-\tfrac12\res_{x,0}(P\log P_2)]\, dx.
\tag3.15$$
\endproclaim

\medskip
Next, we study the second term $-L(P,Q^N_\lambda )$ in (3.2). In the
localization we are considering, 
the operator is built up of s.g.o.s as follows:$$
\aligned
L(P,Q^N_\lambda )&=G^+(P)G^-(Q^N_{\lambda
}),\text{ where}\\
G^+(P)&=r^+Pe^-J,\quad  G^-(Q^N_\lambda
)=Jr^-Q^N_\lambda e^+,\; \text{ with }J\:u(x',x_n)\mapsto
u(x',-x_n);
\endaligned\tag3.16$$
cf.\ \cite{G1, (2.6.5)ff.}. (Here $r^\pm$ denotes
restriction from $\Bbb R^n$ to $\Bbb R^n_\pm$ and $e^{\pm}$
denotes extension by 0 from $\Bbb R^n_\pm$ to $\Bbb R^n$,
respectively.) 
This term contributes in a 
different way than $(PQ_\lambda ^N)_+$. 
It was shown in \cite{GSc1, Sect.\ 3} that the trace expansion of 
$-G^+(P)G^-(Q^N_{\lambda })$ has no term of the form $c(-\lambda
)^{-N}\log (-\lambda )$ and that the coefficient of $(-\lambda
)^{-N}$ is local. The arguments there were based on Laguerre
expansions, which give quite complicated formulas (see e.g.\ Lemmas
A.1 and A.2 in the Appendix), so it is preferable to find a more 
robust method giving
a simpler information on the coefficient of $(-\lambda )^{-N}$.
We shall here use the strategy of \cite{G4, Sect.\ 4} recalled in
Section 2, invoking
the {\it regularity} of parameter-dependent symbols introduced in \cite{G1}.
In the notation of \cite{G1}, the parameter is $\mu =|\lambda
|^\frac12$, where $\lambda $ runs on a ray in $\Bbb C\setminus \rp$,
$-\lambda =\mu ^2e^{i\theta }$ (for some $\theta \in \,]-\pi ,\pi
[\,$).

\proclaim{Lemma 3.4} Consider $G^+(P)G^-(Q_\lambda ^N)$; the order
$\sigma $ of $P$ is integer. 

$1^\circ$ For each $N\ge 1$, the s.g.o.\ $-G^+(P)G^-(Q^N_{\lambda })$ has
order $\sigma -2N$, class $0$ and regularity $\sigma $. Moreover, it can be
written as $$
\aligned
-G^+(P)G^-(Q_{\lambda })&=(-\lambda )^{-1}G^+(P)P_1G^-(Q_{\lambda
})\;\text{ if }N=1,\\
-G^+(P)G^-(Q^N_{\lambda })&=(-\lambda )^{-N}G^+(P)P_1G^{-,N}_\lambda \;\text{
in general},
\endaligned\tag3.17
$$
where $G^+(P)P_1G^{-,N}_\lambda $ has order $\sigma $, class
$0$ and regularity $\sigma +\frac12$. ($G^{-,N}_\lambda $ is described
in {\rm (3.23)} below.)

$2^\circ$ The normal trace $$
\Cal S^{(N)}_\lambda =\tr_n(-G^+(P)G^-(Q^N_{\lambda }))=\tfrac{\partial  _\lambda
^{N-1}}{(N-1)!}\Cal
S^{(1)}_\lambda 
\tag3.18
$$
is a $\psi $do on $\Bbb R^{n-1}$
 of order $\sigma -2N$ and regularity
$\sigma -\frac14$. It can also be
written as $$
\Cal S^{(N)}_\lambda =(-\lambda )^{-N}\Cal S^{(N)\prime}_\lambda
,\quad \Cal S^{(N)\prime}_\lambda =
\tr_n[G^+(P)P_1G^{-,N}_\lambda ],\tag3.19
$$
where $\Cal S^{(N)\prime}_\lambda $ has order $\sigma $ and regularity
$\sigma +\frac14$. 

$3^\circ$ 
The symbol $\frak s^{(N)}(x',\xi ',\lambda )$ of $\Cal
S^{(N)}_\lambda $, with the expansion in (quasi-)homogeneous terms $\frak
s^{(N)}(x',\xi 
',\lambda )\sim \sum_{j\in \Bbb N}\frak s^{(N)}_{\sigma -2N-j}(x',\xi
',\lambda )$ (homogeneous in $(\xi ',\mu )$ on each ray 
$ \lambda =-\mu ^2e^{i\theta
}$, $\mu >0$), satisfies:
$$
\aligned
|\partial _{x',\xi '}^{\beta ,\alpha } [\frak s^{(N)}(x',\xi ',\lambda
)-\sum_{j<J}\frak s^{(N)}_{\sigma   -2N-j}(x',\xi ',\lambda )]|
&\leg \ang{\xi '}^{\sigma   -\frac14 -|\alpha |-J}\ang{\xi ',\mu
}^{-2N+\frac14},\\
|\partial _{x',\xi '}^{\beta ,\alpha } [\frak s^{(N)}(x',\xi ',\lambda
)-\sum_{j<J}\frak s^{(N)}_{\sigma   -2N-j}(x',\xi ',\lambda )]|
&\leg \ang{\xi '}^{\sigma   +\frac14-|\alpha |-J}\ang{\xi ',\mu
}^{-\frac14}\mu ^{-2N},
\endaligned
\tag3.20
$$
on the rays in $\Bbb C\setminus\rp$, for all $\alpha 
,\beta ,J$. 

\endproclaim 

\demo{Proof}
Consider first the case $N=1$. Since $Q_\lambda $ is strongly polyhomogeneous
(the strictly homogeneous symbol is smooth in $\xi $ and $\lambda $ for $|\xi |+|\lambda
|\ne 0$), it is of regularity $+\infty $, and so is the s.g.o.\
$G^-(Q_\lambda )$, of class 0. $G^+(P)$ is $\lambda $-independent of order
$\sigma $ and class 0, hence has regularity $\sigma $ by \cite{G1,
(2.3.54)}, so by the composition rules (cf.\ e.g.\ \cite{G1, (2.6.5)
10$^\circ$}),  the composed operator has order $\sigma -2$, class 0
and regularity $\sigma $. For $N>1$ we have similarly, since $Q_\lambda ^N$ is of order
$-2N$ and regularity $+\infty $, that $G^+(P)G^-(Q^N_{\lambda })$ has
order $\sigma -2N$, class 0 and regularity $\sigma $. This shows the
first statement in $1^\circ$.

For the second statement in $1^\circ$, we use that $$\aligned
(P_1-\lambda )^{-1}&=
(-\lambda )^{-1}-(-\lambda )^{-1}P_1(P_1-\lambda )^{-1}
 \\
(P_1-\lambda )^{-N}&=
(-\lambda
 )^{-N}-(-\lambda )^{-N}\sum_{0\le k< N}c_k(P_1(P_1-\lambda
)^{-1})^{N-k}\\
&=
(-\lambda )^{-N}-(-\lambda )^{-N}P_1\sum_{1\le j\le N}c_{N-j}P_1^{j-1}(P_1-\lambda )^{-j}
\endaligned\tag3.21$$ 
as in (3.10); here $G^-(\lambda ^{-N})=0$. Since $P_1$ is a
differential operator, $G^-(P_1Q_\lambda
)=P_1G^-(Q_\lambda )$, and we find in the case $N=1$:
 $$
-G^+(P)G^-(Q_{\lambda })=(-\lambda )^{-1}G^+(P)P_1G^-(Q_{\lambda }),\tag3.22
$$
where $G^+(P)P_1$ is a $\lambda $-independent s.g.o.\ of order $\sigma +2$ and
class 2, hence has regularity $\sigma +\frac12 $ by \cite{G1,
(2.3.55)}. Then the composed operator $G^+(P)P_1G^-(Q_{\lambda })$ is
of order $\sigma $, class 0 and regularity $\sigma +\frac12$. 
For general $N$ we use that
the last sum in (3.21) is of order $-2$ and regularity $+\infty $. Here
(3.17) holds with$$
G^{-,N}_\lambda =G^-(\sum_{1\le j\le N}c_{N-j}P_1^{j-1}(P_1-\lambda
)^{-j}).\tag 3.23
$$

The statements in $2^\circ$ now follow by use of \cite{G4, Lemma 3.4},
which shows a loss of regularity $\frac 14$ in general when $\tr_n$
is applied.

The information in $3^\circ$ follows directly from the definition of the class
of symbols of the stated regularity,  when $\sigma <0$. When $\sigma
\ge 0$, the regularity information itself gives weaker estimates when
$\sigma -|\alpha |-J\ge 0$; here we use the device from the proof 
of \cite{G4, Prop.\ 4.3}: Compose $\Cal S^{(N)}_\lambda $ to the
left with $\Lambda ^\varrho \Lambda ^{-\varrho }$, $\Lambda
=\operatorname{OP}(\ang{\xi '})$, with $\varrho \ge \sigma +1$. Taking
$\Lambda ^{-\varrho }$ together with $G^+(P)$ one finds that
$\Lambda ^{-\varrho }\Cal S^{(N)}_\lambda $ satisfies the regularity
statements with $\sigma $ replaced by $\sigma -\varrho $, hence the
estimates with the same replacement, and the desired estimates follow
after composition with $\Lambda ^\varrho $.
\qed 
\enddemo

The lemma makes it possible to use Theorems 2.1 and 2.2 as in Section
2.

\proclaim{Theorem 3.5} Consider $\Cal S^{(N)}_\lambda $ defined in
Lemma {\rm 3.4}. When $N>(\sigma +n-1)/2 $, the kernel of $\Cal
S^{(N)}_\lambda $ on the diagonal has an expansion$$
K(\Cal S^{(N)}_\lambda ,x',x')=\sum_{0\le l\le \sigma +n-1}\tilde
{\frak s}^{(N)}_l(x')(-\lambda )
^{\frac{\sigma +n-1 -l}2 -N}+O(\lambda ^{-N-\frac1{8}}),
\tag3.24
$$
where in particular the coefficient of $(-\lambda ) ^{-N}$ is$$
\tilde {\frak s}^{(N)}_{\sigma +n-1 }(x')=
\int_{\Bbb R^{n-1}}{\frak s}^{(N)h}_{ -2N-n+1}(x',\xi
',-1)\,\d\xi'=
\int_{\Bbb R^{n-1}}{\frak s}^{(1)h}_{ -2-n+1}(x',\xi ',-1)\,\d\xi'.\tag3.25
$$

 The ``log-transform'' $S$ with
symbol $s(x',\xi ')$ deduced for $|\xi '|\ge 1$ from the symbol \linebreak
$\frak s^{(1)}(x',\xi 
',\lambda )$ of $ {\Cal S}^{(1)}_\lambda $ by $$
s(x',\xi ')=\tfrac i{2\pi }\int_{\Cal C}\log\lambda \, \frak
s^{(1)}(x',\xi ',\lambda )\, d\lambda 
,\quad s_{\sigma -j}(x',\xi ')=\tfrac i{2\pi }\int_{\Cal C}\log\lambda \, \frak
s^{(1)}_{\sigma -2-j}(x',\xi ',\lambda )\, d\lambda 
\tag3.26
$$
(with a curve $\Cal C$ in $\Bbb C\setminus\crm$ around $[1,\infty
[\,$)
is a classical $\psi $do of order $\sigma $ such that$$
\aligned
\tilde{\frak s}^{(N)}_{\sigma +n-1
}(x')&=-\tfrac12\int_{|\xi '|=1}s_{1-n}(x',\xi ')\,\d S(\xi '),\\
\tr\tilde{\frak s}^{(N)}_{\sigma +n-1
}(x')&=-\tfrac12\operatorname{res}_{x'}S.
\endaligned\tag 3.27
$$
\endproclaim 

\demo{Proof} 
By Lemma 3.4, the operator family $\Cal
S^{(N)}_\lambda $ defined there satisfies the hypotheses for Theorems
2.1 and 2.2  with $m=2$ and $\delta =\frac14$, $V=\Bbb
C\setminus \crp$. For $|\xi '|\ge 1$, the symbols are holomorphic in
$\lambda \in\Bbb C\setminus [1,\infty [\,$. 
\qed
\enddemo 

Let us consider the role of $S$ more closely. Since the coefficient
we are studying is local, it makes no difference if $S$ is modified
for $|\xi '|< 1$. If we replace $P_1$ by $P_1'=\OP([\xi ]^2)$ (cf.\
(1.14)), its resolvent is defined for $\lambda \in \Bbb C\setminus
[\frac12,\infty [\,$, and $\log P'_1$ is well-defined as $\OP(\log
[\xi ]^2)$. 
A simple calculation using \cite{G1, (2.6.44)} shows that
$G^-(Q_\lambda )=$\linebreak $G^-((P_1-\lambda 
)^{-1})$ has symbol-kernel and symbol $$
\aligned
\tilde g^-(q )&=\frac1{2\kappa }e^{-\kappa (x_n+y_n)},\quad \kappa =(|\xi '|^2-\lambda )^{\frac12};\\
g^-(q )&=\frac 1{2\kappa (\kappa +i\xi _n)(\kappa -i\eta
_n)};
\endaligned
\tag3.28 
$$ 
when $|\xi '|\ge 1$ this also holds when $P_1$ is replaced by $P'_1$.

 The symbol terms in $S$ are for $|\xi '|\ge 1$
equal to $$
\aligned
s_{\sigma -j}&(x',\xi ')=\tfrac i{2\pi }\int_{\Cal C}\log\lambda \, \frak
s^{(1)}_{\sigma -2-j}(x',\xi ',\lambda )\, d\lambda \\
&=\tfrac {-i}{2\pi }\int_{\Cal C}\log\lambda \,
\int_{\Bbb R^2} g^+(p)_{\sigma -j}(x',\xi ',\xi _n,\eta _n)\tfrac
1{2\kappa (\kappa +i\eta _n)(\kappa -i\xi 
_n)}\,\d\xi _n\d\eta _nd\lambda \\
&=\tfrac {-i}{2\pi }\int_{\Cal C}\log\lambda \,
\int_{x_n,y_n> 0}\tilde g^+(p)_{\sigma -j}(x',x_n,y_n,\xi ')\tfrac
1{2\kappa 
}e^{-\kappa (x_n+y_n) }\,dx_ndy_nd\lambda .
\endaligned\tag3.29
 $$
Since $\tilde g^+(p)_{\sigma -j}$ is bounded, and $\|\tfrac
1{2\kappa 
}e^{-\kappa (x_n+y_n) }\|_{L_{1,x_n,y_n}({\Bbb R}^2_{++})}\leg \kappa
^{-3}$, hence is $O(\lambda ^{-\frac32})$, the log-integral can be moved inside the
$(x_n,y_n)$-integral. Then we can use, as shown in \cite{GG, Example
2.8}, that in fact 
$$
\tfrac {i}{2\pi }\int_{\Cal C}\log\lambda \,
\tfrac
1{2\kappa 
}e^{-\kappa (x_n+y_n) }\,d\lambda =-\tfrac1{x_n+y_n}e^{-|\xi
'|(x_n+y_n)},
\text{ equal to }\tilde g^-(\log P'_1)\text{
for }|\xi '|\ge 1,
\tag3.30
$$
where $\tilde g^-(\log P'_1)$ is the symbol-kernel of
the generalized s.g.o.\ $G^-(\log P'_1)$. Note that it is integrable in
$(x_n,y_n)\in \rpp $ when $\xi '\ne 0$.

The singularity of the symbol-kernel at $x_n=y_n=0$ is typical for
$G^{\pm}(\log P_2)$ for a general elliptic differential operator $P_2$, as
well as for  the s.g.o.-like part of $\log 
P_{2,T}$ for a realization of $P_2$ defined by a differential elliptic
boundary condition $Tu=0$, cf.\
\cite{GG}.

Now we get, for $|\xi '|\ge 1$,
$$
\multline
s_{\sigma -j}(x',\xi ')=\int_{x_n,y_n>0}\tilde
g^+(p)_{\sigma -j}(x',x_n,y_n,\xi ')
\tfrac
1{x_n+y_n}e^{-|\xi '|(x_n+y_n)}\,dx_ndy_n\\
=-\tr_n(\tilde
g^+(p)_{\sigma -j}\circ_n \tilde g^-(\log P'_1))=
-\operatorname{symb}_{\sigma -j}\tr_n(G^+(P)G^-(\log P'_1));
\endmultline\tag3.31
$$
in short, when we recall (3.16),$$
S\sim-\tr_n(G^+(P)G^-(\log P'_1))= -\tr_n L(P,\log P'_1),\tag3.32
$$
modulo a smoothing operator. With this,
the formula for the coefficient of $(-\lambda )^{-N}$ becomes:
$$
\tilde{\frak s}^{(N)}_{\sigma +n-1
}(x')=\tfrac12 \int_{|\xi
'|=1}\operatorname{symb}_{1-n}(\tr_n L(P,\log P'_1))(x',\xi
')\,\d S(\xi '),
\tag3.33
$$
i.e., for the fiber traces:
 $$
\tr \tilde{\frak s}^{(N)}_{\sigma +n-1
}(x')=-\tfrac12 \operatorname{res}_{x'}S=\tfrac12 \operatorname{res}_{x'}(\tr_nL(P,\log P'_1)).\tag3.34
$$

We conclude:

\proclaim{Theorem 3.6} The operator family $-L(P,(-\Delta -\lambda )^{-N})$ in {\rm (3.2)} on $\Bbb
R^n_+$ has for $N>(\sigma +n)/2$ expansions as in {\rm (1.23)}, with
$\delta =\frac18$, for $\lambda \to\infty $ on rays in $\Bbb
C\setminus \Bbb R_+$. Here $a'_0=0$, and
$$
l_0(-L(P,(-\Delta -\lambda )^{-1}))=
\tfrac12\int_{\Bbb R^{n-1}}\res_{x'}(\tr_n L(P,\log P'_1))\,dx'.  
\tag3.35$$
\endproclaim 

\demo{Proof} This is found by integrating the fiber trace of the
expansion in (3.24) with respect to 
$x'\in\Bbb R^{n-1}$, replacing $l$ by $j-1$, and using the
intepretation of $S$ that we have just accounted for.
\qed
\enddemo

Next, consider the manifold situation, choosing $P_1$ as indicated
after Proposition 3.2. Then, summing over the coordinate 
patches, we get
$l_0(-L(P,(P_1-\lambda )^{-1}))$ 
as a sum of contributions of the form $$
\varphi 
\int \tfrac12
\res_{x'}\tr_n L(P,\log P_1)\,dx'\psi ,\tag3.36$$  
with cutoff functions $\varphi ,\psi $; the integral vanishes if the
local coordinate patch does not meet the boundary.

For a general auxiliary operator $P_2$, we have of course$$
l_0(L(P,(P_2-\lambda )^{-1}))=l_0(L(P,(P_1-\lambda
)^{-1}))+[l_0(L(P,(P_2-\lambda )^{-1}))-l_0(L(P,(P_1-\lambda )^{-1}))].
$$
Here the expression in [\dots] has been treated before, namely in the
course of the proof of \cite{G4, Th.\ 3.6, Rem.\ 3.12}. Rather than
taking up details from that long proof, we simply note that in view of
those results and Theorem 3.3,
$$
\aligned
&l_0(-L(P,(P_2-\lambda )^{-1}))-l_0(-L(P,(P_1-\lambda )^{-1}))\\
&=
C_0(P_+,P_{2,+})-l_0((P(P_2-\lambda )^{-1})_+)
-C_0(P_+,P_{1,+})+l_0((P(P_1-\lambda )^{-1})_+)\\
&=-\tfrac12\res(P_+(\log P_2-\log P_1)_+)+\tfrac12\res((P(\log
P_2-\log P_1))_+)\\
&=\tfrac12\res(L(P,\log P_2-\log P_1)).
\endaligned\tag3.37
$$ The residues here are covered by \cite{FGLS} since $\log P_2-\log
P_1$ is classical having the transmission property, and the expression is
calculated in local coordinates as
$$
\tfrac12\int_{X'}
\res_{x'}(\tr_n L(P,\log P_2-\log P_1))\,dx' .\tag3.38
$$
Combining this with what we found for $l_0(-L(P,(P_1-\lambda )^{-1}))$,
we can conclude:

\proclaim{Theorem 3.7} Assumptions as in Theorem {\rm 3.3}. The
operator family $-L(P,(P_2 -\lambda )^{-N})$  has for $N>(\sigma
+n)/2$ trace
expansions as in {\rm (1.23)}, with 
$\delta =\frac18$, for $\lambda \to\infty $ on rays in $V$. Here $a'_0=0$, and
$$
l_0(-L(P,(P_2 -\lambda )^{-1}))=
\tfrac12\int_{X'}\res_{x'}(\tr_n L(P,\log P_2))\,dx',  
\tag3.39$$
calculated in local coordinates.
\endproclaim

Observe that $\log P_2$ in local coordinates  can be
written as a sum of $\log P'_1$ and a classical $\psi $do 
having the transmission property. Then $G^-(\log P_2)=G^-(\log P'_1)+G^{-,0}$ where $G^{-,0}$ is a
standard singular Green operator (a similar fact is observed 
 in \cite{GG, Prop.\
2.9}). Thus $\frac12\res_{x'}\tr_n L(P,\log
P_2)=\frac12\res_{x'}\tr_n G^+(P)G^-(\log
P_2)$ is a sum of a term with the special function in (3.30) and a term
covered by the residue definition of \cite{FGLS}. 

For general $P_2$, it is therefore a minor extension of the
definitions in \cite{FGLS} to {\it define}
$$
\res L(P,\log P_2)=-2l_0(L(P,(P_2-\lambda )^{-1})).\tag3.40
$$
This number is defined directly in terms of the manifold situation, so
we have ``for free'' that it is independent of local coordinates, 
although we also  have the
description in local coordinates (3.39).

\subhead 4. The constant coming from $GQ_{\lambda ,+}$
\endsubhead

We now study the coefficient of $(-\lambda )^{-N}$ in the expansion of
$GQ^N_{\lambda ,+}= \tfrac{\partial  _\lambda ^{N-1}}{(N-1)!}GQ_{\lambda
,+}$. The order $\sigma $ can here be any real number.
Again we begin with the special localized situation explained in the
beginning of Section 3. 

This is the hardest term to treat, since its
contribution is nonlocal and the homogeneity properties of the
symbols of $G$ and $Q_\lambda $ do not play together in an easy way.  
We here take recourse to one more trace concept, that of ${\Bbb
N}\times{\Bbb N}$-matrices, in the Laguerre expansion representation of the
operator $G$.  

Recall (e.g.\ from \cite{G1}) that $g$ has a rapidly convergent Laguerre expansion 
$$g(x',\xi
',\xi _n,\eta _n)=\sum_{l,m\in\Bbb N}c_{lm}(x',\xi ')\hat\varphi
_l([\xi '],\xi _n)\bar{\hat\varphi }_m([\xi '],\eta _n),\tag4.1$$ 
with
$c_{lm}(x',\xi ')$ polyhomogeneous of order $\sigma $, and Fourier
transformed Laguerre-type functions $$
\hat\varphi
_l([\xi '],\xi _n)=(2[\xi '])^{\frac12}\frac{([\xi ']-i\xi
_n)^l}{([\xi ']+i\xi _n)^{l+1}},\tag4.2$$
$[\xi
']$ defined as in (1.14). In view of the orthonormality of the 
system $\hat \varphi _l$, $l\in\Bbb N$, $$
(\tr_n g)(x',\xi ')=\sum_{l\in\Bbb N}c_{ll}(x',\xi ')\tag4.3
$$
(as used also e.g.\ in \cite{GSc1, (5.17)}), so $\res G$ is defined
from the diagonal terms alone, by integration in $x'$ of $$
\res_{x'}\tr_nG=\int_{|\xi '|=1}\tr \sum_{l\in \Bbb N}c_{ll}(x',\xi ')\,\d S(\xi').\tag4.4
$$
Not just the diagonal terms but also the off-diagonal terms will contribute to
$l_0(G Q_{\lambda ,+})$, so let us define$$\aligned
g_{\operatorname{diag}}&=\sum_{l\in\Bbb N}c_{ll}(x',\xi ')\hat\varphi
_l([\xi '],\xi _n)\bar{\hat\varphi }_l([\xi '],\eta _n),\\ 
g_{\operatorname{off}}&=\sum_{l,m\in\Bbb N, l\ne m}c_{lm}(x',\xi ')\hat\varphi
_l([\xi '],\xi _n)\bar{\hat\varphi }_m([\xi '],\eta _n)=g-g_{\operatorname{diag}},\\ 
\endaligned\tag4.5$$
denoting the corresponding singular Green operators $G_{\operatorname{diag}}$ resp.\  $G_{\operatorname{off}}$.

It is seen from (4.1) that
$$
|g(x',\xi ',\xi _n,\xi _n)|\leg [\xi ']^{\sigma +1}([\xi ']^2+\xi
_n^2)^{-1}, \tag4.6
$$
and there are analogous estimates for homogeneous terms and remainders. One
can furthermore conclude from this that the strictly homogeneous version
$g^h_{\sigma -j}$ (extending $g_{\sigma -j}$ by homogeneity into
$\{1>|\xi '|>0\}$) satisfies$$
|g^h_{\sigma -j}(x',\xi ',\xi _n,\xi _n)|\leg |\xi '|^{\sigma -j+1}(|\xi '|^2+\xi
_n^2)^{-1}. \tag4.7
$$

Let us denote $G(-\Delta -\lambda )^{-N}_+=G'$,
it has symbol $$
g'(x',\xi ',\xi _n,\eta _n,\lambda )=g\circ_n (p_1-\lambda )^{-N}_+=h^-_{-1,\eta _n}[g(x',\xi
',\xi _n,\eta _n)(|(\xi ',\eta _n)
|^2-\lambda)^{-N}],\tag4.8
$$ 
and symbol-kernel $
\tilde g'(x',x_n,y_n,\xi ',\lambda )
=r^+_{x_n}r^+_{y_n}\Cal F^{-1}_{\xi _n\to x_n}\overline{\Cal
F}^{-1}_{\eta  _n\to y_n}g'(x',\xi ',\xi _n,\eta _n,\lambda )$. The
normal trace $\tr_nG'$ is a $\psi $do on $\Bbb R^{n-1}$ with symbol
$$
\aligned
\tr_ng'(x',\xi ',\lambda )&=\int^+g'(x',\xi ',\xi _n,\xi _n,\lambda
)\,\d\xi _n \\
&=\int^+\bigl(h^-_{-1,\eta _n}[g(x',\xi
',\xi _n,\eta _n)(|(\xi ',\eta _n)
|^2-\lambda)^{-N}]\bigr)|_{\eta _n=\xi _n}\d \xi  _n
\endaligned\tag4.9$$
so the kernel of the operator $\tr_nG'$ is, for $x'=y'$, 
$$
\aligned
K(\tr_n G',& x',x',\lambda )\\
&= 
\int_{\Bbb R^{n-1}}\int^+\bigl(h^-_{-1,\eta _n}[g(x',\xi ',\xi
_n,\eta _n)(|(\xi ',\eta 
_n)|^2-\lambda)^{-N}]\bigr)|_{\eta _n=\xi _n}\,\d\xi _n\d\xi '.
\endaligned\tag4.10
$$
We can here make an important simplification, as in \cite{GSc1,
(3.9)ff.}:

\proclaim{Lemma 4.1} Let $r(x',\xi )$ be a $\psi $do
symbol in the calculus, of normal order $ 0$.
For the compositions on the boundary symbol level
$$
g\circ_n r_+=h^-_{-1,\eta _n}[g(x',\xi
',\xi _n,\eta _n)r(x',\xi ',\eta _n)],\quad
r_+\circ_n g=h^+_{\xi _n}[r(x',\xi ',\xi _n)g(x',\xi
',\xi _n,\eta _n)],
$$
one has that
$$
\aligned
\int^+\bigl(h^+_{\eta _n}[g(x',\xi ',\xi _n,\eta _n)r(x',\xi ',\eta 
_n)]\bigr)|_{\eta _n=\xi _n}\,\d\xi _n&=0,\\
\int^+\bigl(h^-_{-1,\xi _n }[r(x',\xi ',\xi _n)g(x',\xi ',\xi _n,\eta
_n)]\bigr)|_{\eta _n=\xi _n}\,\d\xi _n&=0;
\endaligned\tag4.11
$$
and hence
$$\aligned
\tr_n(g\circ _nr_+)&=\int^+\bigl(h^-_{-1,\eta _n}[g(x',\xi
',\xi _n,\eta _n)r(x',\xi ',\eta _n)]\bigr)|_{\eta _n=\xi _n}\d \xi  _n\\
&=\int g(x',\xi ',\xi _n,\xi _n)r(x',\xi ',\xi _n)
\,\d\xi _n;\\
\tr_n(r_+\circ _n g)&=\int^+\bigl(h^+_{\xi  _n}[r(x,\xi ',\xi _n) g(x',\xi
',\xi _n,\eta _n)]\bigr)|_{\eta _n=\xi _n}\d \xi  _n\\
&=\int r(x',\xi ',\xi _n)g(x',\xi ',\xi _n,\xi _n)
\,\d\xi _n.
\endaligned\tag4.12
$$

\endproclaim

\demo{Proof} 
Assume first that $r$ is of normal order $ -1$, hence is $O(\ang{\xi
_n}^{-1})$. 

We have for each term in the Laguerre expansion (4.1),
denoting $[\xi ']=\varrho $:
$$\aligned
&\int^+\bigl(h^+_{\eta _n}[c_{lm}(x',\xi ')\hat\varphi
_l([\xi '],\xi _n)\bar{\hat\varphi }_m([\xi '],\eta _n)r(x',\xi ',\eta
_n)
]\bigr)|_{\eta _n=\xi _n}\,\d\xi _n\\
&=c_{lm}\int^+\hat\varphi _l(\varrho ,\xi _n)
\bigl(h^+_{\eta _n}[
\bar{\hat\varphi }_m(\varrho ,\eta _n)
r(x',\xi ',\eta  _n)
]\bigr)
|_{\eta _n=\xi _n}\,\d\xi _n,\\
&=c_{lm}\int^+\hat\varphi _l(\varrho ,\xi _n)
h^+_{\xi  _n}[ 
\bar{\hat\varphi }_m(\varrho ,\xi _n)
r(x',\xi ',\xi   _n)]
\,\d\xi _n.
\endaligned$$
Here the integrand is holomorphic in $\xi _n$ on $\Bbb
C_-$ and $C^\infty $ on $\overline{\Bbb C}_-$, and is $O(\xi _n^{-2})$ for $|\xi _n|\to\infty $ in
$\overline{\Bbb C}_-$, 
so that the integration curve can be transformed to a closed curve in
$\Bbb C_-$ where it gives 0 (one may check this with \cite{G1,
(2.2.42)ff.}). Summing over $l$ and $m$, we find the first line in
(4.11), and the first statement in (4.12) is an immediate consequence,
since $h^+f(\xi _n)+h^-_{-1}f(\xi _n)=f(\xi _n)$ for $f\in \Cal
H_{-1}$.

Similarly, for the second line in (4.11), each Laguerre term is a
standard integral of an $L_2$ function that is seen to be 0 by changing
the integration curve to a closed curve in ${\Bbb C}_+$, and the
second statement in (4.12) follows.

Now if $r$ is of normal order $0$, it can be written as a sum 
$
r(x',\xi )=r_0(x',\xi ')+r'(x',\xi )
$,
where $r_0$ is independent of $\xi _n$ (and polynomial in $\xi '$),
and $r'$ is of normal order $-1$, thanks to the transmission
condition. Here the $\circ_n$ composition with $r_0$ is simply a
multiplication that goes outside the $\tr_n$-integral, and the
preceding considerations apply to $r'$.
\qed\enddemo

Thus $h^-_{-1,\eta _n}$ can be omitted in (4.10),
and we arrive at the more convenient formula
$$
\aligned
K(\tr_n G', x',x',\lambda )
&= 
\int_{\Bbb R^{n-1}}\int g(x',\xi ',\xi _n,\xi _n)
(|(\xi ',\xi  
_n)|^2-\lambda)^{-N}\,\d\xi _n\d\xi '\\
&=\int_{\Bbb R^{n}}g(x',\xi ',\xi _n,\xi _n)(|(\xi ',\xi  
_n)|^2-\lambda)^{-N}\,\d\xi .
\endaligned\tag4.13
$$

For this, the analysis of the asymptotic expansion in $-\lambda $
will to some extent be modeled after the proof of Proposition 3.1; but
it presents additional difficulties since the terms in $g$ are given
as homogeneous for $|\xi '|\ge 1$ only, and polar coordinates in
$\xi $ are not very helpful. 
In relation to the expansion $g\sim \sum_{j\in\Bbb N}g_{\sigma -j}$ in
homogeneous terms (with $g_{\sigma -j}(x',\xi ',\xi _n,\eta _n)$ homogeneous in $(\xi
',\xi _n,\eta _n)$ {\it of degree} $\sigma -j-1$ for $|\xi '|\ge 1$,
hence of {\it order} $\sigma -j$, in the notation of  \cite{G4}), 
we set $g_{1-n}=0$ if $\sigma +n-1\notin\Bbb N$, and we define $
g_{<1-n}=g-\sum_{\sigma -j\ge 1-n}g_{\sigma -j}$.

The terms with $\sigma -j\ne 1-n$ are relatively easy to handle, and
will be treated first:

\proclaim{Lemma 4.2} For $g_{\sigma -j}(x',\xi ',\xi _n, \eta _n)$
with $\sigma -j>1-n$, one has:$$
\aligned
\int_{\Bbb R^n}g_{\sigma -j}&(x',\xi ',\xi _n,\xi _n)
(|\xi |^2-\lambda )^{-N}\,\d\xi \\
&=(-\lambda
)^{\frac{\sigma +n-1-j}2-N}\int_{\Bbb R^n}g^h_{\sigma -j}(x',\xi ',\xi
_n,\xi _n)(|\xi
|^2+1)^{-N}\,\d\xi \\
&\quad+(-\lambda )^{-N}\slint (\tr_ng_{\sigma -j})(x',\xi
')\,\d\xi '+O(\lambda ^{-N-\frac12+\varepsilon }),
\endaligned\tag4.14
$$
for $\lambda \to\infty $ on rays in $\Bbb C\setminus\Bbb R_+$. For 
$g_{<1-n}(x',\xi ',\xi _n, \eta _n)$ one has:$$
\int_{\Bbb R^n}g_{<1-n}(x',\xi ',\xi _n,\xi _n)
(|\xi |^2-\lambda )^{-N}\,\d\xi 
=(-\lambda )^{-N}\slint (\tr_ng_{<1-n})(x',\xi
')\,\d\xi '+O(\lambda ^{-N-\frac12+\varepsilon }).\tag4.15
$$
\endproclaim

\demo{Proof} 
Let $\sigma -j>1-n$. It follows from (4.7) that $g^h_{\sigma -j}$
is integrable on cylinders $\{\xi \in\Bbb R^n\mid |\xi '|\le
a\}$. Multiplication by $(|\xi |^2-\lambda
)^{-N}$ makes it integrable over $\Bbb R^n$ 
when $\lambda \notin \crp$. 
Write
$$
\aligned
\int_{\Bbb R^n}g_{\sigma -j}&(x',\xi ',\xi _n,\xi _n)(|\xi |^2-\lambda )^{-N}\,\d\xi \\
&=\int_{\Bbb R^n}g^h_{\sigma -j}(x',\xi ',\xi _n,\xi _n)(|\xi |^2-\lambda )^{-N}\,\d\xi \\
&\quad+\int_{|\xi '|\le 1}(g_{\sigma -j}(x',\xi ',\xi _n,\xi _n)-g^h_{\sigma -j}(x',\xi ',\xi
_n,\xi _n))(|\xi |^2-\lambda )^{-N}\,\d\xi . 
\endaligned\tag4.16
$$

For the first term we have by homogeneity for $\lambda $ on rays in
$\Bbb C\setminus{\overline{\Bbb R}_+}$, writing $\lambda =-|\lambda |e^{i\theta }$,
$|\theta |<\pi $, and replacing $\xi $ by $|\lambda |^{\frac12}\eta $:
$$
\multline
\int_{\Bbb R^n}g^h_{\sigma -j}(x',\xi ',\xi _n,\xi _n)(|\xi
|^2+|\lambda |e^{i\theta } 
)^{-N}\,\d\xi 
\\=|\lambda |^{\frac{\sigma -j+n-1}2-N}\int_{\Bbb R^n}g^h_{\sigma
-j}(x',\eta ',\eta _n,\eta _n )(|\eta
|^2+e^{i\theta }
)^{-N}
\,\d\eta .\endmultline\tag4.17 
$$ 
This equals the first term in the right hand side of (4.14) if $\theta
=0$, and the identity extends analytically to general $\lambda $ (as
in \cite{GS1, Lemma 2.3}).

Using (3.10),
we find that
$$
\multline
\int_{|\xi '|\le 1}(g_{\sigma -j}(x',\xi ', \xi _n,\xi _n)-g^h_{\sigma -j}(x',\xi ', \xi _n,\xi _n))(|\xi
|^2-\lambda )^{-N}\,\d\xi \\
=(-\lambda )^{-N}\int_{|\xi '|\le 1}(g_{\sigma -j}(x',\xi ', \xi _n,\xi _n)-g^h_{\sigma
-j}(x',\xi ', \xi _n,\xi _n))
\,\d\xi +O(\lambda ^{-N-\frac12+\varepsilon }),\endmultline\tag4.18
$$
since $(g_{\sigma -j}-g^h_{\sigma -j})\ang{\xi _n}^{1-2\varepsilon }$
is integrable on the set where $|\xi '|\le 1$ (cf.\ (4.6), (4.7)).
Here,  in view of (2.2), and (3.6) applied in dimension $n-1$,$$
\multline
\int_{\xi \in\Bbb R^n,|\xi '|\le 1}(g_{\sigma -j}(x',\xi ', \xi _n,\xi _n)-g^h_{\sigma
-j}(x',\xi ', \xi _n,\xi _n))
\,\d\xi \\
=\int_{\xi '\in\Bbb R^{n-1},|\xi '|\le 1}((\tr_ng_{\sigma -j})(x',\xi
')-(\tr_ng^h_{\sigma 
-j})(x',\xi '))
\,\d\xi '=\slint (\tr_ng_{\sigma -j})(x',\xi ')\,\d\xi ';
\endmultline
\tag4.19
$$
this shows (4.14). 

For (4.15), we use that $g_{<1-n}$ is integrable in $\xi $ and
remains so after multiplication by a small power of $|\xi |$; then
(3.10) can be used to show that
$$
\int g_{<1-n}(x',\xi ', \xi _n,\xi _n)(|\xi
|^2-\lambda )^{-N}\,\d\xi =(-\lambda )^{-N}\int g_{<1-n}(x',\xi ', \xi _n,\xi _n)
\,\d\xi +O(\lambda ^{-N-\frac12+\varepsilon }),
$$
where
$$
\int g_{<1-n}(x',\xi ', \xi _n,\xi _n)
\,\d\xi 
=\slint (\tr_ng_{<1-n})(x',\xi ')\,\d\xi ',
$$
in view of (2.2), and (3.6) in dimension $n-1$.
\qed
\enddemo

If $\sigma \notin \Bbb Z$ or $\sigma <1-n$, this ends the analysis,
since there is no 
term in $g$ of order $1-n$. Since$$
\slint (\tr_n g)(x',\xi ')\,\d\xi '=\sum_{\sigma -j>1-n}\slint
(\tr_ng_{\sigma -j})(x',\xi ')\,\d\xi '+\slint (\tr_ng_{<1-n})(x',\xi
')\,\d\xi '
$$ 
in this case, we thus find:

\proclaim{Theorem 4.3} When $\sigma -n+1\notin \Bbb N$, the kernel of
$\tr_n (GQ^N_{\lambda ,+})$ on the 
diagonal $\{x'=y'\}$ has for $N>(\sigma +n)/2$ an expansion:$$
\multline
K(\tr_n(GQ^N_{\lambda ,+}),x',x')=\sum_{0\le j<\sigma +n-1}(-\lambda
)^{\frac{\sigma +n-1-j}2-N}b_j(x')\\
+(-\lambda
)^{-N}\slint (\tr_n g)(x',\xi ')\,\d \xi  + O(\lambda
^{-N-\frac12+\varepsilon })\endmultline
\tag4.20
$$ 
any $\varepsilon >0$, for $\lambda \to\infty $ on rays in $\Bbb
C\setminus \Bbb R_+$, where$$
b_j(x')=\int_{\Bbb R^{n}} g^h_{\sigma -j}(x',\xi ',\xi _n,\xi _n)(|\xi
|^2+1)^{-N}\,\d \xi.
\tag4.21$$
\endproclaim

When $\sigma $ is an integer $\ge 1-n$, we must include a study of the
contribution from $g_{1-n}$. Here
the Laguerre expansion (4.1) comes into the picture, since we need to
treat $G_{\operatorname{diag}}$ and  $G_{\operatorname{off}}$ (cf.\ (4.5))
by different methods. The symbols of order $1-n$ are:
$$
\gathered
g_{\operatorname{diag},1-n}=\sum_{l\in{\Bbb N}}g_{ll,1-n},\quad
g_{\operatorname{off},1-n}=\sum_{l,m\in{\Bbb N}, l\ne
m}g_{lm,1-n},\text{ where}\\
g_{lm,1-n}(x',\xi
',\xi _n,\eta _n)=c_{lm,1-n}(x',\xi ')\hat\varphi
_l([\xi '],\xi _n)\bar{\hat\varphi }_m([\xi '],\eta _n).
\endgathered\tag4.22
$$
We denote $\OP(c_{lm}(x',\xi '))=C_{lm}$.

\proclaim{Proposition 4.4} 
The $g_{ll,1-n}$ satisfy for $N\ge 1$:
$$
\aligned
\int_{\Bbb R^n}&\tr g_{ll,1-n}(x',\xi ', \xi _n,\xi _n)(|\xi |^2-\lambda
)^{-N}\,\d\xi  \\
&=\tfrac{\partial  _\lambda ^{N-1}}{(N-1)!}\int_{\Bbb R^n}\tr g_{ll,1-n}(x',\xi ', \xi _n,\xi _n)(|\xi |^2-\lambda
)^{-1}\,\d\xi  \\
&=\tfrac{\partial  _\lambda ^{N-1}}{(N-1)!}\bigl[\tfrac12\res_{x'}C_{ll}(-\lambda
)^{-1}\log(-\lambda )+\slint \tr c_{ll,1-n}(x',\xi
')\,\d\xi '(-\lambda )^{-1}\\
&\quad +\res_{x'}C_{ll}(
-\log 2 \,
(-\lambda 
)^{-1}+O(\lambda
^{-\frac32+\varepsilon }))\bigr],\\
&=\tfrac12 \res _{x'}C_{ll}(-\lambda
)^{-N}\log(-\lambda )+\slint \tr
c_{ll,1-n}(x',\xi ')\,\d\xi '(-\lambda)^{-N}\\
&\quad +\res_{x'}C_{ll}\bigl(-(\log 2 +\tfrac12 \alpha _N)(-\lambda
)^{-N}+O(\lambda
^{-N-\frac12+\varepsilon })\bigr) ,
\endaligned \tag4.23$$
for $\lambda \to\infty $ on rays in $\Bbb C\setminus\Bbb R_+$.

It follows by summation in $l$ that
$$
\aligned
\int_{\Bbb R^n}&\tr g_{\operatorname{diag},1-n}(x',\xi ', \xi _n,\xi _n)(|\xi |^2-\lambda
)^{-N}\,\d\xi  \\
&=\tfrac{\partial  _\lambda ^{N-1}}{(N-1)!}[\tfrac12\res_{x'}(\tr_n G)(-\lambda
)^{-1}\log(-\lambda )\\
&\quad+\bigl(\slint\tr
\tr_n g_{1-n}(x',\xi ')\,\d\xi '-\log 2 \,\res_{x'}(\tr_n G)
\bigr)(-\lambda 
)^{-1}]+O(\lambda
^{-N-\frac12+\varepsilon }),\\
&=\tfrac12 \res _{x'}(\tr_nG)(-\lambda
)^{-N}\log(-\lambda )\\
&\quad+\bigl(\slint\tr
\tr_n g_{1-n}(x',\xi ')\,\d\xi '-(\log 2 +\tfrac12 \alpha
_N)\res_{x'}(\tr_n G) \bigr)(-\lambda
)^{-N}\\
&\quad +O(\lambda
^{-N-\frac12+\varepsilon }).
\endaligned \tag4.24$$
\endproclaim

\demo{Proof}
Since the integrand is in $L_1(\Bbb R^n)$ for all $N\ge 1$, we can
insert the formula\linebreak
$(|\xi |^2-\lambda )^{-N}=\tfrac{\partial  _\lambda
^{N-1}}{(N-1)!}(|\xi |^2-\lambda )^{-1}$ and pull the differential
operator $\tfrac{\partial  _\lambda ^{N-1}}{(N-1)!}$ outside the
integral sign.
Then we can do the main work in the case $N=1$.
Write
$$
\aligned
\int_{\Bbb R^n}\tr g_{ll,1-n}&(x',\xi ', \xi _n,\xi _n)(|\xi |^2-\lambda
)^{-N}\,\d\xi=I_1+I_2,  
\\I_1&=\int_{|\xi '|\ge 1 }\tr g^h_{ll,1-n}(x',\xi ', \xi _n,\xi _n)(|\xi
|^2-\lambda )^{-N}\,\d\xi,\\
I_2&=\int_{|\xi '|\le 1}\tr g_{ll,1-n}(x',\xi ', \xi _n,\xi _n)
(|\xi |^2-\lambda )^{-N}\,\d\xi .
\endaligned\tag4.25$$
The last term is treated as in the preceding proof. Using (3.10), and the
orthonormality of the Laguerre functions,
we find that
$$\aligned
I_2&=\int_{|\xi '|\le 1}\tr g_{ll,1-n}(x',\xi ', \xi _n,\xi _n)(|\xi |^2-\lambda
)^{-N}\,\d\xi \\
&=(-\lambda )^{-N}\int_{|\xi '|\le 1}\tr g_{ll,1-n}(x',\xi ', \xi _n,\xi _n)
\,\d\xi +O(\lambda ^{-N-\frac12+\varepsilon })\\
&=(-\lambda )^{-N}\slint \tr c_{ll,1-n}(x',\xi ')
\,\d\xi '+O(\lambda ^{-N-\frac12+\varepsilon }),
 \endaligned\tag4.26$$
in view of (2.2), and (3.6) in dimension $n-1$.

Now consider $I_1$. It suffices to take $\lambda
\in\Bbb R_-$; here we write $-\lambda =\mu ^2$, $\mu >0$.
The integrand is homogeneous in $(\xi ',\xi _n)$ but the integration
is over $\{\xi \mid |\xi '|\ge 1\}$; it is here that the Laguerre
expansion helps in the calculations.
Write in general 
$$
I_{lm}=\int_{|\xi '|\ge 1}\tr c_{lm,1-n}(x',\xi ')\hat\varphi
_l([\xi '],\xi _n)\bar{\hat\varphi }_m([\xi '],\xi  _n)(|\xi |^2-\lambda )^{-N}\,\d\xi .\tag4.27
$$
Denote $r=|\xi '|$, $ \kappa =(|\xi '|^2-\lambda )^{\frac12}=(r^2+\mu
^2)^{\frac12}$. Then by Lemma A.2,
$$
I_1=I_{ll}=\deln \res_{x'}C_{ll}
\int_{r\ge 1}\frac1{r\kappa (r +
\kappa )}
\,dr.\tag4.28
$$
 The integral is analyzed by use of Lemma A.3. Here
(A.4) implies:$$
\aligned
&I_{ll}
=\res_{x'}C_{ll}\tfrac{\partial _\lambda ^{N-1}}{(N-1)!}
[\mu ^{-2}(-\log 2+\log(\sqrt{1+\mu ^2}\,+1))]
\\
&=\res_{x'}C_{ll}\tfrac{\partial _\lambda ^{N-1}}{(N-1)!}
[(-\lambda )^{-1}(-\log 2+\log((-\lambda
)^{\frac12})+\log(\sqrt{1+(-\lambda )^{-1}}\,+(-\lambda )^{-\frac12}))]
\\
&=\res_{x'}C_{ll}\tfrac{\partial _\lambda ^{N-1}}{(N-1)!}
[\tfrac12(-\lambda )^{-1}\log(-\lambda
)-\log 2(-\lambda )^{-1}+O(\lambda ^{-\frac32})],
\endaligned\tag4.29$$
 for $|\lambda |\to\infty $.

When we add $I_1$ and $I_2$, we
obtain the statement of 
(4.23) before the last identity. In the last identity, the term
$-\frac12 \alpha _N\operatorname{res}_{x'}(C_{ll})$ (cf.\ (1.5)) 
comes from differentiating the log-term,
and it is checked from (4.29) that the error has the asserted order.

Now summation with respect to $l$ gives the second statement in the
proposition, in view of (3.16)ff.
\qed
\enddemo

The coefficient $\log 2$ in (4.23), (4.24) may seem a little odd, but
fortunately, it will disappear again when the terms are analyzed more
and set in relation to compositions with $\log P'_1$.
First we make the appropriate
analysis of $g_{\operatorname{off},1-n}$.
 An analysis similar to the above of the $I_{lm}$ with $l\ne m$ 
can of course be carried
out, but as seen from (A.3) and (A.5),  this gives
coefficients and remainders depending on $l$ and $m$, 
leading to a less transparent
summation result.

We shall use instead that this part
in fact has good enough regularity 
properties to define a completely local contribution, where the methods
of Sections 2 can be applied. It was observed already in
\cite{GSc1}
that this part gives a local coefficient and no logarithmic term.

\proclaim{Proposition 4.5} 
Define for $N\ge 1$,$$
\Cal S^{(N)}_{\operatorname{off},\lambda
}
=\tr_n(G_{\operatorname{off}}\deln(P_1-\lambda )_+^{-1}),\text{ with
symbol }\frak s^{(N)}_{\operatorname{off}};\tag4.30
$$
it is of order $\sigma -2N$ and regularity $\sigma -\frac14$, and it
also equals $(-\lambda)
^{-N}$ times a symbol
of order $\sigma $ and regularity $\sigma
+\frac14$, and satisfies estimates like {\rm (3.20)}. Then we can define
the ``log-transform'' $S_{\operatorname{off}}$ with symbol
$$
s_{\operatorname{off}}(x',\xi ')=\tfrac i{2\pi }\int_{\Cal C}\log\lambda \, \frak
s_{\operatorname{off}}(x',\xi ',\lambda )\, d\lambda \tag4.31
$$
(for $|\xi '|\ge 1$).

The fiber trace of the kernel of 
$\Cal S^{(N)}_{\operatorname{off},\lambda }
$ has an expansion on the diagonal 
$$
\tr K(\Cal S^{(N)}_{\operatorname{off},\lambda } ,x',x')=\sum_{0\le l\le
\sigma +n-1}\tr \tilde
{\frak s}^{(N)}_{\operatorname{off},l}(x')(-\lambda )
^{\frac{\sigma +n-1 -l}2 -N}+O(\lambda ^{-N-\frac1{8}}),
\tag4.32
$$
where
$$
\tr\tilde{\frak s}^{(N)}_{\operatorname{off},\sigma +n-1
}(x')=-\tfrac12\operatorname{res}_{x'}S_{\operatorname{off}}.
\tag4.33
$$
\endproclaim 

\demo{Proof}
For the symbol $\frak s_{\operatorname{off}}(x',\xi ',\lambda )$ of $\Cal S_{\operatorname{off},\lambda }$,
the immediate rules give that $\frak s^{(N)}_{\operatorname{off}}(x',\xi ',\lambda )=\deln \frak s_{\operatorname{off}}(x',\xi ',\lambda )$ is
of order $\sigma -2N$ and regularity $\sigma -\frac14$. But
it can be rewritten using (3.21) (we let $|\xi '|\ge 1$ and as usual
write $|\xi '|=r$ for short): 
$$\aligned
\frak s_{\operatorname{off}}(x',\xi ',\lambda )&=\sum_{l\ne m}c_{lm}(x',\xi
')\tr_n[\hat\varphi _l(r,\xi 
_n)\overline{\hat\varphi }_m(r,\eta _n)\circ_n (p_1-\lambda )_+^{-1}]\\
&=\sum_{l\ne m}c_{lm}(x',\xi ')\tr_n[\hat\varphi _l(r,\xi
_n)\overline{\hat\varphi }_m(r,\eta _n)\circ_n (-\lambda
^{-1}+\lambda ^{-1}p_1(p_1-\lambda )^{-1})_+]\\
&=\sum_{l\ne m}c_{lm}(x',\xi ')\tr_n[\hat\varphi _l(r,\xi
_n)\overline{\hat\varphi }_m(r,\eta _n)\circ_n ( p_1\lambda
^{-1}(p_1-\lambda )^{-1})_+], 
\endaligned$$
where we used that $\tr_n(\hat\varphi _l(r,\xi
_n)\overline{\hat\varphi }_m(r,\eta _n))=(\hat\varphi _l,{\hat\varphi
}_m)=0$ for $l\ne m$. Similarly to Lemma 
3.2, we see from this that $\Cal S_{\operatorname{off},\lambda }$ can also be viewed as $\lambda
^{-1}$ times an operator of order $\sigma $ and regularity $\sigma
+\frac14$, and likewise  $(-\lambda )^{N}\deln \Cal S_{\operatorname{off},\lambda }$ is
of order $\sigma $ and regularity $\sigma
+\frac14$; moreover, estimates like (3.20) hold. 
Then Theorems 2.1 and 2.2 apply, and the result of
Theorem 3.5 extends to this case, when we define $s_{\operatorname{off}}$ from $\frak s_{\operatorname{off}}$
by (4.31).
\qed 
\enddemo 

These two propositions can be taken together to give a formula for the
term $l_0(G Q_{\lambda ,+})$. However, to find a more transparent formulation, we shall
analyze the ingredients some more.

Insertion of $\frak s_{\operatorname{off}}$ in (4.31) gives that (for
$|\xi  '|\ge 1$)
$$
s_{\operatorname{off}}(x',\xi ')=\tfrac i{2\pi }\int_{\Cal C}\log\lambda \, 
\sum_{l\ne m}c_{lm}(x',\xi ')\tr_n[\hat\varphi _l(r,\xi
_n)\overline{\hat\varphi }_m(r,\eta _n)\circ_n (p_1 \lambda
^{-1}(p_1-\lambda )^{-1})_+]\,d\lambda . 
$$
One would like to pass the integration in $\lambda $ inside $\tr_n$,
so that $\log p_1$ would appear in the formula, but the complex rule
for calculation of $\tr_n(g\circ_n r_+)$ (cf.\ (4.9), (4.12)) requires
a symbol $r$ satisfying the transmission condition, which $\log p_1$
does not (in its principal part). However, inspired from (4.12) 
we can {\it define} an
extension of the normal trace to non-standard symbols:

\proclaim{Definition 4.6} For functions $r(x',\xi )$ and $g(x',\xi
,\eta _n)$, we set
$$
\aligned
\tr'_n(g\circ_n r_+)&=\int_{{\Bbb R}}g(x',\xi ',\xi _n,\xi _n)r(x',\xi
',\xi _n)\,\d\xi _n,\\
\tr'_n(r_+\circ_n g)&=\int_{\Bbb R} r(x',\xi ',\xi _n)
g(x',\xi ',\xi _n,\xi _n )
\,\d\xi _n ,
\endaligned\tag4.34
$$
whenever the integral converges. The notation will also be used for
the associated operators. 
\endproclaim

This allows us to write
$$\aligned
s_{\operatorname{off}}(x',\xi ')&=\tfrac i{2\pi }\int_{\Cal C}\log\lambda \, 
\sum_{l\ne m}c_{lm}(x',\xi ')\tr'_n[\hat\varphi _l(r,\xi
_n)\overline{\hat\varphi }_m(r,\eta _n)\circ_n (p_1 \lambda
^{-1}(p_1-\lambda )^{-1})_+]\,d\lambda\\
&=\tfrac i{2\pi }\int_{\Cal C}\log\lambda \, 
\sum_{l\ne m}c_{lm}(x',\xi ')\int_{{\Bbb R}}\hat\varphi _l(r,\xi
_n)\overline{\hat\varphi }_m(r,\xi  _n) p_1(\xi ) \lambda
^{-1}(p_1(\xi )-\lambda )^{-1}\,\d\xi _nd\lambda \\ 
&=
\sum_{l\ne m}c_{lm}(x',\xi ')\int_{{\Bbb R}}\hat\varphi _l(r,\xi
_n)\overline{\hat\varphi }_m(r,\xi  _n) p_1(\xi ) \tfrac i{2\pi }\int_{\Cal C}\log\lambda \, \lambda
^{-1}(p_1(\xi )-\lambda )^{-1}\,d\lambda  \d\xi _n\\
&=
\tr'_n (g_{\operatorname{off}}\circ_n  (p_1 
 (p_1^{-1}\log p_1))_+)
=
\tr'_n (g_{\operatorname{off}}\circ_n  (\log p_1)_+).
\endaligned\tag4.35$$
In this sense, we have that
$$S_{\operatorname{off}}\sim\tr'_n(G_{\operatorname{off}}(\log P'_1)_+)
\tag4.36$$
(modulo smoothing operators).

We can now extend this point of view to cases where
$G_{\operatorname{off}}$ is replaced by $G_{\operatorname{diag}}$ or
$G$ itself.
{\it Define:}
$$
\aligned
S_{\operatorname{diag}}&\sim\tr'_n(G_{\operatorname{diag}}(\log
P'_1)_+),\\
S&\sim\tr'_n(G(\log
P'_1)_+),
\endaligned
\tag4.37$$
with symbols given for $|\xi '|\ge 1$ by:
$$
\aligned
s_{\operatorname{diag}}(x',\xi ')&= \sum_{l\in\Bbb N}c_{ll}(x',\xi ')\tr'_n[\hat\varphi _l(r,\xi
_n)\overline{\hat\varphi }_l(r,\eta _n)\circ_n  (\log
p_1)_+],\\
s(x',\xi ')&= \sum_{l,m\in\Bbb N}c_{lm}(x',\xi ')\tr'_n[\hat\varphi _l(r,\xi
_n)\overline{\hat\varphi }_m(r,\eta _n)\circ_n  (\log
p_1)_+]\\
&= \tr'_n[g(x',\xi ',\xi _n,\eta  _n)\circ_n  (\log
p_1)_+].
\endaligned\tag4.38
$$

Further calculations show:

\proclaim{Lemma 4.7} The symbol $s_{\operatorname{diag}}$ of
$S_{\operatorname{diag}}$ satisfies$$
s_{\operatorname{diag}}(x',\xi ')\sim (\tr_n g)(x',\xi ')(2\log 2+
2\log[\xi ']) .\tag4.39  
$$
It is log-polyhomogeneous, and in particular,$$
\res_{x',0}(S_{\operatorname{diag}})=\res_{x',0}(\tr'_n(G_{\operatorname{diag}}(\log
P'_1)_+))=2\log 2\,
\res_{x'}(\tr_nG).\tag4.40 
$$
\endproclaim 
\demo{Proof} For $r=|\xi '|\ge 1$, write $$
\aligned
s_{\operatorname{diag}}(x',\xi ')&=\sum_{l\in\Bbb N}s_{ll}(x',\xi
'),\text{ with}\\
s_{ll}(x',\xi ')&=c_{ll}(x',\xi ')\tr'_n[\hat\varphi _l(r,\xi
_n)\overline{\hat\varphi }_l(r,\eta _n)\circ_n  (\log
p_1)]\\
&=c_{ll}(x',\xi ')\int_{{\Bbb R}}\hat\varphi _l(r,\xi
_n)\overline{\hat\varphi }_l(r,\xi  _n)  \log (r^2+\xi _n^2)\,\d\xi _n.
\endaligned$$
Here we have
$$
\multline
\int_{{\Bbb R}}\hat\varphi _l(r,\xi
_n)\overline{\hat\varphi }_l(r,\xi  _n)  \log (r^2+\xi _n^2)\,\d\xi
_n=
\int _{\Bbb R}2r \frac{(r - it)^{l}}{(r + it)^{l+1}}\frac{(r +
it)^{l}}{(r -
it)^{l+1}}\log(r^2+t^2)\,\d t\\
=\frac{2r}{2\pi }\int_{\Bbb R}\frac {\log(r^2+t^2)}{r^2+t^2}\,dt=\frac1{\pi} \int
_{\Bbb R}\frac{\log(1+s^2)+\log r^2}{1+s^2}\,ds
=2\log 2+\log r^2
\endmultline$$
(see Lemma A.4 in the Appendix). It follows that $$
s_{\operatorname{diag}}(x',\xi ')=\sum_{l\in\Bbb N}c_{ll}(x',\xi
')(2\log 2+2\log |\xi '|),\tag4.41
$$
for $|\xi '|\ge 1$, cf.\ also (4.3). This shows (4.39),
and (4.40) follows in view of (4.4).\qed
\enddemo

It follows from this and (4.35) that $S=\tr'_n(G(\log
P'_1)_+)$  is likewise log-polyhomogeneous.

Then we can finally conclude, in the special localized situation:

\proclaim{Theorem 4.8} The operator family $G(-\Delta -\lambda )^{-N}_+$ on $\Bbb
R^n_+$ has for $N>(\sigma +n)/2$, expansions as in {\rm (1.23)}, with
$\delta =\frac12+\varepsilon $, for $\lambda \to\infty $ on rays in $\Bbb
C\setminus \Bbb R_+$. 
Here $a'_0=\tfrac12\res(G)=\tfrac12\int_{{\Bbb R}^{n-1}}\res_{x'}\tr_nG\,dx'$ and$$
l_0(G(-\Delta -\lambda )_+^{-1})=
\int_{{\Bbb R}^{n-1}}\bigl(\TR_{x'}\tr_n G-\tfrac12\res_{x',0}\tr'_n(G(\log P'_1)_+)\bigr)\,dx'.    
\tag4.42$$
\endproclaim 
 
\demo{Proof} 
Collecting the results of Lemma 4.2, Theorem 4.3 for $\sigma +n-1\notin{\Bbb N}$ (where the relevant
residues vanish), and
Propositions 4.4 and 4.5 for $\sigma +n-1\in{\Bbb N}$, we get
the expansion with $a'_0$ as stated and,  $$
\multline
 l_0(G(-\Delta -\lambda )_+^{-1})=\int_{{\Bbb R}^{n-1}}\bigl(\TR_{x'}\tr_n G-\log 2\,\res_{x'}\tr_nG
-\tfrac12\res_{x'}S_{\operatorname{off}}\bigr)\,dx'\\    
=\int_{{\Bbb R}^{n-1}}\bigl(\TR_{x'}\tr_n G-\log 2\,\res_{x'}\tr_nG
-\tfrac12\res_{x'}\tr'_n(G_{\operatorname{off}}(\log P'_1)_+)\bigr)\,dx',    
\endmultline\tag4.43$$ 
in view of (4.36).
Now since $G=G_{\operatorname{diag}}+G_{\operatorname{off}}$ with 
$\tr'_n(G(\log P'_1)_+)$ log-polyhomogeneous,$$
\aligned
\res_{x'}\tr'_n(G_{\operatorname{off}}(\log P'_1)_+)
&=\res_{x',0}\tr'_n(G(\log P'_1)_+)-\res_{x',0}\tr'_n(G_{\operatorname{diag}}(\log P'_1)_+)\\
&=\res_{x',0}\tr'_n(G(\log P'_1
)_+)-2\log 2\,\res_{x'}\tr_n
G, 
\endaligned
$$
by Lemma 4.7. When this is inserted in (4.43),
we find (4.42).\qed 
\enddemo 

This gives the formula for the contribution to the basic zeta
coefficient in a reasonably natural form. Now the global formula with
$P_1$, as well as the formula with $P_1$
replaced by general $P_2$, follow as in
Section 3:

\proclaim{Theorem 4.9} For a singular Green operator 
$G$ of order $\sigma \in\Bbb R$ and class $0$, together with an
elliptic differential operator $P_2$ as in Theorem {\rm 3.3},
$G(P_2-\lambda )^{-N}_+$ has 
expansions as in {\rm (1.23)} for $N>(\sigma +n)/2$, with
$\delta =\frac12+\varepsilon $, for $\lambda \to\infty $ on rays in
$V$.
Here $a'_0=\tfrac12\res(G)$
and
 $$
l_0(G(P_2-\lambda )^{-1}_+) =\int_{X '}\bigl(\TR_{x'}
\tr_n G-\tfrac12\res_{x',0}\tr'_n(G(\log P_2)_+)\bigr)\,dx',
\tag4.44
$$
calculated in local coordinates.
\endproclaim 

It is worth keeping in mind here that $\tr'_n(G(\log P_2)_+))$ in
local coordinates is the sum of a usual normal trace, namely that of $G(\log
P_2-\log P_1)_+$, and a generalized normal trace $\tr'_n$ in one
special case, namely for $G$
composed with the logarithm of the Laplacian.
\medskip

\subhead 5. The constants coming from $P_+G_\lambda $ and $GG_\lambda $\endsubhead

The appropriate expansion of $\Tr(P_++G)G^{(N)}_\lambda $ was shown in
Section 2, and we shall just discuss the constants arising there
more in depth.
In this analysis $m$ can be any positive integer, and we as usual take
$N>(\sigma +n)/m$. 

Recall that the operator
$U$ introduced in Corollary 2.5 arose from taking the normal trace of
$BG^{(N)}_\lambda $ and integrating together with $\log\lambda $; this
gave a $\psi $do on $X'$. We shall set this in relation to the
composition of $B$ with the s.g.o.-like
part $G_1^{\log}$ of $\log P_{1,T}$:
$$
G_1^{\log}=\lim _{s\searrow 0}\tfrac i{2\pi }\int_{\Cal C}\lambda
^{-s}\log \lambda \,G_\lambda \, d \lambda =\tfrac i{2\pi }\int_{\Cal
C}\lambda ^{-1}\log \lambda \,P_1G_\lambda \, d \lambda \tag5.1 
$$ 
(cf.\ (2.13)); it is described in more detail in \cite{GG}. 
The analysis is 
different for the two terms coming from $P_+$ resp.\ $G$, so let us
write
$$
U=U_P+U_G,\tag 5.2
$$
where $U_P$ is defined from $B=P_+$ and $U_G$ is defined from $B=G$.
We have the trace expansions from Corollary 2.5:
$$\aligned
\Tr(P_+G^{(N)}_\lambda )&=\sum_{1\le j\le \sigma +n}a^{(N)}_j(-\lambda )
^{\frac{\sigma +n-j}m -N}+O(\lambda ^{-N-\frac1{4 m}})\text{ with }
\\
a^{(N)}_{\sigma +n}&=-\tfrac1m \operatorname{res}U_P \;\big(
=l_0(P_+G_\lambda )\big); 
\\
\Tr(GG^{(N)}_\lambda )&=\sum_{1\le j\le \sigma +n}b^{(N)}_j(-\lambda )
^{\frac{\sigma +n-j}m -N}+O(\lambda ^{-N-\frac1{4 m}(+\varepsilon
)})\text{ with }
\\
 b^{(N)}_{\sigma +n}&=-\tfrac1m \operatorname{res}U_G \;\big( =l_0(GG_\lambda )\big).
\endaligned\tag5.3$$

The natural thing to do would be to interchange the $\lambda $-integral
and the application of $\tr_n$. This works well for the part with
$G$. Here we can proceed as in the discussion of $S$ in Section 3,
by appealing to the estimates of the symbol-kernel of $G_\lambda $ shown by
Seeley \cite{S} (and recalled in \cite{GG, Th.\ 2.4}):
The
symbol-kernel has in local coordinates an expansion
in quasi-homogeneous terms $\tilde g\sim \sum_{j\ge 0}\tilde g_{-m-j}$, 
 satisfying estimates on the rays in $V$, with $ \kappa =| \xi '|+|\lambda |^{\frac1m}$:
$$
|D_{x'}^\beta D_{\xi'}^\alpha x_n^kD_{x_n}^{k'}
y_n^lD_{y_n}^{l'}D_\lambda ^p
\tilde g_{-m-j}|\leg \kappa ^{1-m-|\alpha |-k+k'-l+l'-j-mp}e^{-c\kappa
(x_n+y_x)}\tag5.4 
$$
for all indices, when $ \kappa \ge \varepsilon $.
These estimates allow carrying the $\lambda $-integration inside
$\tr_n$:

\proclaim{Theorem 5.1} For $U_G$ in {\rm (5.2)}, we have that
$U_G=\tr_n(G_1G^{\log})$, and the coefficient of $(-\lambda )^{-N}$ in
the trace expansion of $GG^{(N)}_\lambda $ is:
$$
l_0(GG_\lambda )=-\tfrac1m \res U_G= -\tfrac1m \operatorname{res} \tr_n (GG_1^{\log}).\tag5.5
$$
\endproclaim 

\demo{Proof}
Consider the $j$'th symbol term $u_{\sigma -j}(x',\xi ')$ in $U_G$, cf.\
(2.22). It
is for $|\xi '|\ge 1$ a linear combination of terms of the form $$
\aligned
&\tfrac {i}{2\pi }\int_{\Cal C}\log\lambda \,
\int_{x_n,y_n> 0} \partial_{\xi '}^\alpha  \tilde g_{\sigma -k}(x',x_n,y_n,\xi ')
\partial_{x'}^\alpha \tilde g_{-l-m}(x',y_n,x_n,\xi ',\lambda )
\,dx_ndy_nd\lambda \\
&
=\int_{x_n,y_n> 0} \partial_{\xi '}^\alpha  \tilde g_{\sigma -k}(x',x_n,y_n,\xi ')
\tfrac {i}{2\pi }\int_{\Cal C}\log\lambda \,\partial_{x'}^\alpha \tilde g_{-l-m}(x',y_n,x_n,\xi ',\lambda )
\,d\lambda dx_ndy_n,
\endaligned\tag5.6
 $$
for $|\alpha |+k+l=j$ ($\alpha \in{\Bbb N}^{n-1}, k,l\in {\Bbb N}$); we could interchange the order of integration
since the symbol-kernels from $G$ are bounded and those from $G_\lambda $
have ${L_{1,x_n,y_n}({\Bbb R}^2_{++})}$-norms that are 
$O(\lambda ^{-1-\frac1m})$. The inner log-integrals define the terms
in the symbol-kernel of $G_1^{\log}$ (and derivatives), cf.\ \cite{GG}.
We conclude that $\tr_n(GG_1^{\log})$ is a well-defined classical
$\psi $do of order $\sigma $ which equals $U_G$ modulo a
smoothing operator.\qed
\enddemo

Here it is only a slight extension of the definition of \cite{FGLS} to
define
$$
\res(GG_1^{\log})=-m\,l_0(GG_\lambda ).\tag5.7
$$

For $U_P$, it would be nice to be able to write a similar relation of
$\res U_P$ to the residue of $\tr_n(P_+G_1^{\log})$. But the latter
normal trace rarely exists in the usual sense; it does so when $P$ is
of low order, but otherwise needs interpretations involving e.g.\ the
subtraction of a number of symbol terms which do not
contribute to the residue. This is apparent already in the case $P=I$,
as considered in \cite{GG, Sect.\ 3}, where $\tr_n$ is defined for
$G_1^{\log}$ minus its principal part.

For clarity in the explanation, assume first that 
in the local coordinates in ${\Bbb R}^n_+$, where  we are considering
$P_+G_1^{\log}$, 
the symbol $p$ is
independent of $x_n$. Then, in view of the transmission condition, it can be decomposed into a differential
operator symbol and a symbol of normal order $-1$:
$$
\aligned
p(x',\xi )&=p'(x',\xi )+p''(x',\xi ), \\ 
p'(x',\xi )&=\sum_{|\alpha
|\le \sigma }a_\alpha (x')(\xi ')^{\alpha '}\xi _n^{\alpha _n},\quad
|p''|\leg \ang{\xi _n}^{-1}\ang{\xi '}^{\sigma +1}.
\endaligned\tag5.8
$$
(One can instead let $p''$ be the part of normal order 0.)  For the differential operator part, the composition of
$D_{x_n}^{\alpha _n}$ with $G_1^{\log}$ is simply an application of
$D_{x_n}^{\alpha _n}$ to the symbol-kernel $\tilde g^{\log}$. This
makes the symbol-kernel {\it more} singular at $x_n+y_n=0$, the larger $\alpha _n$
is, cf.\ \cite{GG, (2.23)--(2.24)}; this is also clear from the
example
$\tilde g^{\log}=\frac1{x_n+y_n}e^{-|\xi '|(x_n+y_n)}$ where
$P_1=-\Delta $.  However, the
term that contributes to the residue escapes this effect since the
order is lifted by $\alpha _n$ so that we have to go further down in
the series to find the term whose order equals minus the boundary dimension.
In short, the relevant term in the $\circ_n$-composition of $p'$ and
$\tilde g^{\log}$ is
$$
\sum_{|\alpha |\le \sigma ,|\alpha '|+\alpha _n-j=1-n}
a_\alpha (x')(\xi ')^{\alpha '}D_{x_n}^{\alpha _n}\tilde g^{\log}_{-j}(x',x_n,y_n,\xi '), \tag5.9
$$
which is $O((x_n+y_n)^{-\varepsilon })$ when $n=2$, bounded in
$x_n+y_n$ when $n>2$, and integrable for $x_n+y_n\to \infty $,
according to \cite{GG,(2.26)--(2.29)}. In the full composition there
will moreover be terms coming from derivatives in $\xi '$ and $x'$:
$$
\sum_{|\alpha |\le \sigma ,\beta '\le \alpha ',|\alpha'-\beta ' |+\alpha _n-j=1-n}
c_{\alpha ',\beta '}a_\alpha (x')(\xi ')^{\alpha '-\beta '}D_{x_n}^{\alpha
_n}\partial_{x'}^{\beta '}\tilde g^{\log}_{-j}(x',x_n,y_n,\xi '), \tag
5.10
$$
but we still have estimates as in  \cite{GG,(2.26)--(2.29)}, since
$|\alpha '-\beta '|\ge 0$, so $\tr_n$ {\it of the relevant terms} makes
sense. In this sense we identify $\res U_{P'}$ with $\res \tr_n (P'_+G^{\log}_1)$. 

For the part $p''$ of normal order $-1$ (or 0) we can use in the
$\circ_n$-composition with the symbol-kernel of $G_\lambda $ that   
$\tr_n$ of it can be simplified, by Lemma 4.1. 
Now apply the notation $\tr_n'$ from Definition 4.6.
Considering the terms in 
the full $\circ$-composition of $p''$ with $g(x',\xi ',\xi _n,\eta
_n,\lambda )$, we see that the $\log \lambda $-integration can be
interchanged with taking $\tr'_n$, leading to an identification of
$U_{P''}$ with $\tr'_n(P''_+G_1^{\log})$ 
on the symbol level, so that
$\res U_{P''}$ identifies with $\res\tr'_n(P''_+G_1^{\log})$.

In the case where the symbol of $P$ depends also on $x_n$, one
applies the preceding considerations to the terms in a Taylor
expansion in $x_n$ (a technique used
extensively in \cite{G1, Ch.\ 2.7}), up to a term with a so large
power of $x_n$ that the order is so low that there is no residue
contribution. This just gives some more terms in slightly 
more complicated formulas; the principle is the same as above.

This is as close as we can get to identifying $\res U_P$ with the
residue of a normal trace associated with
$P_+G_1^{\log}$. But it gives sufficient justification for {\it defining}
$$
\res_{x'}(P_+G_1^{\log})=\int_{|\xi '|=1}\tr\,[\tr_n
\operatorname{symb}_{1-n}(P'_+G_1^{\log})+\tr'_n
\operatorname{symb}_{1-n}(P''_+G_1^{\log})]\, \d S(\xi '),
\tag5.11
$$
consistently with $\res_{x'}U_P
$, with reference to the detailed interpretations of the
contributions from $P'$ and $P''$ given
above. Then we can write: 

\proclaim{Theorem 5.2} For $U_P$ in {\rm (5.2)},
the coefficient of $(-\lambda )^{-N}$ in the trace expansion of
$P_+G^{(N)}_\lambda $ is (cf.\ {\rm (5.11))}:
$$
l_0(P_+G_\lambda)=-\tfrac1m \res U_P= -\tfrac1m \int_{X'}\operatorname{res} _{x'}
(P_+G_1^{\log})\, dx',\tag5.12
$$
also denoted $-\frac1m \res(P_+G_1^{\log})$.
\endproclaim

\medskip

\subhead 6. Consequences \endsubhead

We collect the results of Sections 3--5 in the following theorem,
denoting the general auxiliary operator $P_1$:

\proclaim{Theorem 6.1}
Let $B=P_++G$ on $X$, of order $\sigma $, where $P$ is a
$\psi $do satisfying the transmission
condition 
at $X'$ and $G$ is a singular Green operator of class $0$ ($\sigma \in
{\Bbb Z}$ if $P\ne 0$, $\sigma \in{\Bbb R}$ otherwise). Moreover, let $P_1$ be
a second-order elliptic differential operator
provided with a differential boundary condition $Tu=0$,
such that $P_{1,T}$ has $\rmi$ as a spectral cut, and let $V$ be a sector
around $\rmi$ such that the principal symbol
and principal boundary symbol realization have no eigenvalues in $V$.
The operator family $B(P_{1,T} -\lambda )^{-N}=$ \linebreak$(P_++G)(Q^N_{\lambda ,+}+G^{(N)}_\lambda )$ has for $N>(\sigma
+n)/2$ trace expansions as in {\rm (1.23)}, for $\lambda
\to\infty $ on rays in $V$. 
The basic zeta coefficient $C_0(B, P_{1,T})=l_0(B(P_{1,T} -\lambda
)^{-1})$ is a sum of five invariantly defined terms:
$$
\aligned
C_0(B, P_{1,T})&=l_0((PQ_\lambda )_+)-l_0(L(P,Q _\lambda ))+l_0(GQ _{\lambda ,+})\\
&\quad+l_0(P_+G_\lambda )+l_0(GG_\lambda );
\endaligned\tag6.1
$$
here
$$
\aligned
l_0((PQ _\lambda )_+)&=\int_X[\TR _xP-\tfrac12\operatorname{res}_{x,0}(P\log
P_1)]\, dx,\\
-l_0(L(P,Q _\lambda ))&
=\tfrac12\operatorname{res} L(P,\log P_1)=\tfrac12\int_{X'}\res_{x'}\tr_n L(P,\log  P_1))\,dx' ,\\
l_0(GQ _{\lambda ,+})&=\int _{X'}[\TR _{x'}\tr_n
G-\tfrac12\operatorname{res}_{x',0}\tr'_n(G(\log P_1)_+)]\, dx',\\ 
l_0(P_+G_\lambda )&=-\tfrac12\res(P_+G_1^{\log})\\
 =-\tfrac12&\int_{X'}\int_{|\xi '|=1}\tr\,[\tr_n
\operatorname{symb}_{1-n}(P'_+G_1^{\log})+\tr'_n
\operatorname{symb}_{1-n}(P''_+G_1^{\log})]\, \d S(\xi '),\\
l_0(GG_\lambda )&=
-\tfrac12\res(GG_1^{\log})=-\tfrac12\int_{X'}\res_{x'}\tr_n (GG_1^{\log})\,dx' ,\\
\endaligned
\tag6.2$$
where the integrals are calculated in local coordinates ($P$ is
decomposed as in {\rm (5.8)ff}). 
\endproclaim

Note that the first and the third terms in (6.1) are global, and
closely related to the expression (1.15) found in \cite{PS} for the
boundaryless case; the other terms are local.
\medskip

Let us now also show how this can be generalized to
auxiliary operators of even order $m=2k$. The goal is to obtain formulas
as in (6.1)--(6.2) with $\frac12$ replaced by $\frac1m=\frac1{2k}$.

First choose the second-order operator $P_{1,T}$ as a selfadjoint
positive operator (e.g.\ the 
Dirichlet realization  of $P_0+P^*_0+a$, where $P_0$ has principal
symbol $|\xi |^2$ and $a$ is sufficiently large). Then
$(P_{1,T})^k=(P_1^k)_{T'}$ 
with trace operator $T'=\{T,TP_1,\dots,TP_1^{k-1}\}$; this is likewise 
positive selfadjoint, and by definition of $C_0$,
$$
C_0(B,(P_{1,T})^k)=C_0(B,P_{1,T}).\tag6.3
$$
 For the latter, we have the formulas (6.1), (6.2). Since $\log P_1^k=k\log
P_1$, $\frac12\log P_1$ can be replaced by $\frac1m\log P_1^k$ in the
formulas, so the first and third term, the nonlocal terms, are
completely analogous to the case 
$m=2$. The last two terms in (6.1) have already been shown to be of
the desired form with $\frac12$ replaced by $\frac1m$, since Sections
2 and  5 treat arbitrary $m$. For the second term, one can establish
an analysis similar to that in Section 3, using that the symbol-kernel
of $G^-(P^k_1)$ is  of the form $-\frac k{x_n+y_n}e^{-|\xi '|(x_n+y_n)}$
plus a standard singular Green symbol-kernel (for $|\xi '|\ge 1$),
cf.\ \cite{GG, Prop.\ 2.7}.

This shows a generalization of Theorem 6.1 for the special choice
$(P^k_1)_{T'}$. The analysis is now extended to more general auxiliary
operators $P_{2,T''}$ of order $m=2k$ using that
$$
\multline
C_0(B, P_{2,T''})=C_0(B, P_{2,+})+l_0(BG_\lambda )\\
=C_0(B, (P^k_1)_+)+[C_0(B, P_{2,+})-C_0(B, (P^k_1)_+)]
+l_0(BG_\lambda );
\endmultline\tag6.4$$
here Theorem 2.6 shows
formula (1.12) for 
$C_0(B,P_{2,+})-C_0(B, (P^k_1)_+)
$, and the contributions from the singular Green part
$G_\lambda $ of $(P_{2,T''}-\lambda )^{-1}$ are as worked out in
Section 5. We then obtain:

\proclaim{Corollary 6.2} Theorem {\rm 6.1} holds with $P_1$ of even
order $m$, when $\frac12$ in the formulas is replaced by $\frac1m$,
as written in {\rm (1.25)}.
\endproclaim

As a further consequence of these results, 
we can formulate some defect formulas (relative formulas) which  compare two different
choices of auxiliary operator, generalizing (1.12). For one thing, we
can supply (1.12) with a formula where the full auxiliary operators,
not just their $\psi $do parts, enter:

\proclaim{Corollary 6.3} For $B$ together with two auxiliary operators
$P_{1,T_1}$ and $P_{2,T_2}$ as in Corollary {\rm 6.2}, of even orders
$m_1$ resp.\ $m_2$, one has:
$$
\multline
C_0(B,P_{1,T_1})-C_0(B,P_{2,T_1})=
-\operatorname{res}(B(\tfrac1{m_1} \log
P_1-\tfrac1{m_2}\log P_2)_+)\\
+\operatorname{res} (L(P,\tfrac1{m_1}\log
P_1- \tfrac1{m_2}\log P_2))
-\tfrac1{m_1}\res(P_+G_1^{\log})+\tfrac1{m_2}\res(P_+G_2^{\log})\\
-\tfrac1{m_1}\res(GG_1^{\log})+\tfrac1{m_2}\res(GG_2^{\log}),
\endmultline\tag 6.5
$$
with residues defined as above.
\endproclaim  

\demo{Proof} We write $\log P_{i,T_i}=(\log P_i)_++G_i^{\log}$ for
$i=1,2$. Then (6.5) follows simply by applying Corollary 6.2 to the
corresponding zeta coefficients and forming the difference.\qed
\enddemo

The new terms in comparison with (1.12) are local, defined by an
extension of the residue
concept as explained at the end of Section 5.
\medskip

Another defect that we can consider is the difference between the zeta
coefficients arising from different choices of the definition of
$\log$.
It may happen that the spectrum of $P_{1,T}$ (and that of $P_1$
on $\widetilde X$) divides the complex
plane into several sectors with infinitely many eigenvalues, separated
by rays without eigenvalues. Consider two such rays $e^{i\theta }{\Bbb
R}_+$ and $e^{i\varphi }{\Bbb R}_+$, with $\theta <\varphi < \theta 
+2\pi $. We can then define $\log_\theta P_{1}$  and $\log_\theta P_{1,T}$ 
by variants of (1.13),
where the logarithm is taken with cut at $e^{i\theta }{\Bbb R}_+$ and
the integration curve goes around this cut (more details in \cite{GG});
the logarithmic term in the symbol of $\log_\theta P_1$ is still
$m\log[\xi ]$.
The difference between the logarithms is closely connected with a
spectral projection; in fact,
$$
\aligned
\log_\theta P_{1}-\log_\varphi P_{1}&= \tfrac{2\pi }i \Pi _{\theta
,\varphi }(P_{1})\text{ on }\widetilde X,\\
\log_\theta P_{1,T}-\log_\varphi P_{1,T}&= \tfrac{2\pi }i \Pi _{\theta
,\varphi }(P_{1,T})\text{ on }X,
\endaligned
\tag6.6
$$
where $\Pi _{\theta ,\varphi }(P_{1})$  and $\Pi _{\theta ,\varphi }(P_{1,T})$ 
are the sectorial projections for
the sector $\Lambda _{\theta ,\varphi }$, essentially projecting onto
the generalized eigenspace for eigenvalues of $P_1$ resp.\ $P_{1,T}$ in  $\Lambda _{\theta ,\varphi }$
along the complementing eigenspace (see \cite{GG} for the precise
description of such sectorial projections). Here the
$\psi $do part of $\Pi _{\theta ,\varphi }(P_{1,T})$  is $\Pi _{\theta
,\varphi }(P_{1})_+$ satifying the transmission condition since $m$ is
even, and 
the s.g.o.-type part $G_{\theta ,\varphi }$ of $\Pi _{\theta ,\varphi
}(P_{1,T})$ is connected with the s.g.o.-type part of the logs by
$$
G_1^{\log_\theta }-G_1^{\log_\varphi  }=\tfrac {2\pi }i G_{\theta
,\varphi }.\tag 6.7
$$

Since the expressions $\TR_x P$ and $\TR_{x'}\tr_nG$ are independent
of the auxiliary operator, the terms resulting from these are
the same when one calculates the zeta coefficients 
$$
C_{0,\theta }(B, P_{1,T})=l_{0,\theta }(BR_\lambda ),  \text{ resp. }
C_{0,\varphi }(B, P_{1,T})=l_{0, \varphi  }(BR_\lambda ),
\tag6.8
$$
referring to expansions for $\lambda \to \infty $ along $e^{i\theta
}{\Bbb R}_+$ resp.\ $e^{i\varphi }{\Bbb R}_+$. So they disappear when
we take the difference between these zeta coefficients. This implies
another local defect formula:

\proclaim{Corollary 6.4} For the difference between the basic zeta
coefficients  calculated
with respect to two different rays $e^{i\theta }{\Bbb R}_+$ and
$e^{i\varphi  }{\Bbb R}_+$, for $P_1$ of even order $m$, we have:
$$
\multline
C_{0,\theta }(B, P_{1,T})-C_{0,\varphi  }(B, P_{1,T})=
-\tfrac1m\res_+(P(\log_\theta P_1- \log _\varphi P_1))\\
+\tfrac1m\res(L(P,(\log_\theta P_1- \log _\varphi P_1)))
-\tfrac1m \res (G(\log_\theta P_1- \log _\varphi P_1)_+)\\
-\tfrac1m\res(P_+(G_1^{\log_\theta }-G_1^{\log_\varphi }))
-\tfrac1m\res(G(G_1^{\log_\theta }-G_1^{\log_\varphi }));
\endmultline
$$
here the first three terms are residues in the sense of \cite{FGLS},
and the last two are generalizations defined as in {\rm (6.2)}.
In view of {\rm (6.6)--(6.7)}, this can also be written:
$$
\multline
C_{0,\theta }(B, P_{1,T})-C_{0,\varphi  }(B, P_{1,T})=
-\tfrac{2\pi}{im}\res_+(P\,\Pi _{\theta ,\varphi }(P_1))\\
+\tfrac{2\pi}{im}\res(L(P,\Pi _{\theta ,\varphi }(P_1)))
-\tfrac{2\pi}{im} \res (G(\Pi _{\theta ,\varphi }(P_1))_+)
-\tfrac{2\pi}{im}\res((P_++G)G_{\theta ,\varphi }),
\endmultline\tag6.9
$$
in short denoted $-\tfrac{2 \pi
}{im}\res(B\, \Pi _{\theta ,\varphi }(P_{1,T}))$
(the residues being defined from the homogeneous terms of the relevant
dimension as in {\rm (6.2)}). 
When $\Pi _{\theta ,\varphi }(P_{1,T})$
belongs to the standard calculus, this is a residue in the sense of
\cite{FGLS}. 
 
\endproclaim        

Finally we have some observations on the traciality of the new residue
definitions. Let $B$ be {\it of order and class} 0; it is a bounded operator
in $L_2$, and its adjoint is of the same
kind. Then it can be placed to the right of
$R_\lambda ^{N}$ in
all the above calculations, performed in the analogous
way. So $\Tr(R_\lambda ^{N}B)$ has an expansion similar to (1.3) (this
also follows directly from applying (1.3) to $B^*(P_{1,T}^*-\bar\lambda
)^{-N}$ and conjugating). Then since $\Tr(BR_\lambda
^{N})=\Tr(R_\lambda ^{N}B)$, it follows that
$$
l_0(BR_\lambda )=l_0(R_\lambda B),
\text{ in particular }
l_0(BQ_{\lambda ,+})=l_0(Q_{\lambda ,+}B),
l_0(BG_\lambda )=l_0(G_\lambda B).
\tag6.11
$$
This leads to:

\proclaim{Theorem 6.5} When $B=P_++G$ is of order and class
$0$ and $P_{1,T}$ is as in Corollary {\rm 6.2}, then 
$$
 \res(BG_1^{\log}) =\res(G_1^{\log}B),\quad \res([B,\log P_{1,T}])=0.\tag6.12 
$$
Moreover, in the situation of Corollary {\rm 6.4}, 
$$
\aligned
\res(B\,G _{\theta ,\varphi })&=\res(G _{\theta  ,\varphi }\,
B)\\
\res(B\,\Pi _{\theta ,\varphi }(P_{1,T}))&=\res(\Pi _{\theta  ,\varphi }
(P_{1,T})\,B).
\endaligned
\tag6.13
$$
\endproclaim

\demo{Proof} First we note that $\Tr(G^{(N)}_\lambda B)$ has an
expansion similar to (2.25), since the normal trace of
$G^{(N)}_\lambda B$ allows application of Theorems 2.1 and 2.2 as in
the case of $BG^{(N)}_\lambda $. The coefficient of $(-\lambda )^{-N}$
is $l_0(G_\lambda B)=-\tfrac1m \operatorname{res}U'  $, where the
$\psi $do $U'$ on $X'$ is defined similarly as in Corollary 2.5 from taking the normal trace of
$G^{(N)}_\lambda B$ and integrating together with $\log\lambda $.
The constant $\res U'$ is interpreted as
$\res (G_1^{\log }B)$ in the same way as in Section 5; for the part $\res
U'_G$ this goes directly by an interchange of the integration in
$\lambda $ and the application of $\tr_n$, and for the part $\res
U'_P$ it works as in the consideration of the contribution from $P''$,
where the $\lambda $-integration and the application of
$\tr'_n$ to the relevant terms could be interchanged. The first formula in
(6.12) then follows from $l_0(BG_\lambda )=l_0(G_\lambda B)$.

For the second formula in (6.12), we do not have a separate residue
definition for $B\log P_{1,T}$ (neither for $\log P_{1,T}B$), since (6.2)
contains global terms. However, the $\TR$-components of $P$ and $G$ 
will be the same in the formulas for $l_0(BR_\lambda )$ and
$l_0(R_\lambda B)$, so they cancel out in the difference, leaving a
local contribution; it is this that we denote $\res([B,\log
P_{1,T}])$. It vanishes because of the first identity in (6.11).

In the situation of Corollary 6.4,
the preceding results give the commutativity with
$G_1^{\log_\theta }$ or $G_1^{\log_\varphi }$ inserted; this
implies the first line in (6.13) in view of (6.7). Now the second
line follows since the $\psi $do part of  $\Pi _{\theta ,\varphi
}(P_{1,T})$ is a standard $\psi $do in the calculus, for which the
property is known from \cite{FGLS}.
\qed
\enddemo

\subhead Appendix \endsubhead

Some lemmas on integrals arising from Laguerre expansions will be
included here.

\proclaim{Lemma A.1} Set
$$
\gathered
r=|\xi '|,\quad \kappa =(|\xi '|^2-\lambda )^{\frac12}=(r^2+\mu
^2)^{\frac12};\\
s^{\pm}_{l,m}(\xi ',\mu )=\int _\Bbb R\frac 
{(r - i\xi
_n)^{l}}{(r + i\xi
_n)^{l+1}}\frac{(r +
i\xi  
_n)^{m}}{(r -
i\xi  
_n)^{m+1}}\frac1{\kappa
\pm i\xi _n}\,\d\xi _n,
\endgathered
\tag A.1
$$
then one has:
$$\alignedat2
&\text{For $l>m$, }&&s^+_{l,m}=0,\quad
s^-_{l,m}=\frac{(r -\kappa )^{l-m-1}}{(r +\kappa 
)^{l-m+1}}.\\
&\text{For $l=m$, }&&s^+_{l,l}=\frac1{(\kappa +r ) 2r
}= s^-_{l,l}.\\
&\text{For $l<m$, }&&s^+_{l,m}=\frac{(r -\kappa
)^{m-l-1}}{(r +\kappa  
)^{m-l+1}},\quad s^-_{l,m}=0.
\endalignedat\tag A.2$$

\endproclaim

This was shown in \cite{GSc1, Lemma 5.2}.

\proclaim{Lemma A.2} Let $c_{lm,1-n}(x',\xi ')$ be homogeneous of
degree $1-n$ in $\xi '$ for $|\xi '|\ge 1$. Then
$$\aligned
I_{lm}&=\int_{|\xi '|\ge 1}c_{lm,1-n}(x',\xi ')2r\tfrac{(r - i\xi
_n)^{l}}{(r + i\xi
_n)^{l+1}}\tfrac{(r +
i\xi  
_n)^{m}}{(r -
i\xi  
_n)^{m+1}}(r^2+\xi _n^2-\lambda )^{-N}
\,\d\xi \\
&=\int_{r\ge 1}\int_{t\in\Bbb R}2\tfrac{(r - it)^{l-m-1}}{(r + it)^{l-m+1}} 
\tfrac{\partial  _\lambda ^{N-1}}{(N-1)!} \Bigl(\tfrac1{2\kappa
}\Bigl(\tfrac1{\kappa +it}+\tfrac1{\kappa -it}\Bigr)\Bigr)
dr\d t\int_{|\xi '|=1}c_{lm,1-n}(x',\xi
')\d S(\xi ').
\endaligned$$
 Here
$$\aligned
\int_{r\ge 1}\int_{t\in\Bbb R}&\tfrac{(r - it)^{l-m-1}}{(r +
it)^{l-m+1}}
\tfrac{\partial _\lambda ^{N-1}}{(N-1)!} \Bigl(\tfrac1{\kappa
}\Bigl(\tfrac1{\kappa +it}+\tfrac1{\kappa -it}\Bigr)\Bigr)
\,dr\d t\\
&=\tfrac{\partial _\lambda ^{N-1}}{(N-1)!}
\int_{r\ge 1}\int_{t\in\Bbb R}\tfrac{(r - it)^{l-m-1}}{(r +
it)^{l-m+1}}
 \Bigl(\tfrac1{\kappa
}\Bigl(\tfrac1{\kappa +it}+\tfrac1{\kappa -it}\Bigr)\Bigr)
\,dr\d t\\
&=\cases \tfrac{\partial _\lambda ^{N-1}}{(N-1)!}
\int_{r\ge 1}\tfrac{(r -\kappa )^{l-m-1}}{\kappa (r +
\kappa )^{l-m+1}}
 \,dr\text{ if }l> m,\\
\tfrac{\partial _\lambda ^{N-1}}{(N-1)!}
\int_{r\ge 1}\tfrac1{r\kappa (r +
\kappa )}
\,dr\text{ if }l= m,\\
\tfrac{\partial _\lambda ^{N-1}}{(N-1)!}
\int_{r\ge 1}\tfrac{(r -\kappa )^{m-l-1}}{\kappa (r +
\kappa )^{m-l+1}}
 \,dr\text{ if }l< m.\endcases
\endaligned
\tag A.3$$
\endproclaim

\demo{Proof} The first calculation is straightforward, using that  
$c_{lm,1-n}(\xi ')=r^{1-n}c_{lm,1-n}(\xi '/|\xi '|)$. The reductions
of the integrals 
over $r$ and $t$ in (A.3) now follow from (A.2).\qed 
\enddemo

The integrals in (A.3) can be calculated by use of the following
lemma:

\proclaim{Lemma A.3} Let $a>0$. 
One has that $$
\int_1^\infty \tfrac1{ r\sqrt{r^2+a^2}\,(\sqrt{r^2+a^2}\,+r)}\,
dr=a^{-2}(-\log 2+\log(\sqrt{1+a^2}\,+1)),\tag A.4
$$
and, for $j=1,2,\dots$,
$$
\int_1^\infty \tfrac{(\sqrt{r^2+a^2}\,-r)^{j-1}}{
\sqrt{r^2+a^2}\,(\sqrt{r^2+a^2}\,+r)^{j+1}}\, 
dr=(2j)^{-1}a^{-2}(\sqrt{1+a^{-2}}\,-a^{-1})^{2j}.\tag A.5
$$ 
\endproclaim 

\demo{Proof} Letting $s=r/a$, we have:$$\aligned
\int_1^\infty &\tfrac1{ r\sqrt{r^2+a^2}\,(r+\sqrt{r^2+a^2}\,)}\,dr
=a^{-2}\int_{1/a}^\infty \tfrac1{ s\sqrt{s^2+1}\,(\sqrt{s^2+1}\,+s)}\,
ds\\
&=a^{-2}\int_{1/a}^\infty \tfrac{\sqrt{s^2+1}\,-s\,}{ s\sqrt{s^2+1}\,}
\,ds
=a^{-2}\int_{1/a}^\infty \Bigl(\tfrac1s-\tfrac1{ \sqrt{s^2+1}\,}\Bigr)
\,ds\\
&=a^{-2}\bigl[\log s-\log (\sqrt{s^2+1}\,+s)\bigr]_{1/a}^\infty 
=-a^{-2}\bigl[\log (\sqrt{1+s^{-2}}\,+1)\bigr]_{1/a}^\infty \\
&=a^{-2}(-\log 2+\log(\sqrt{1+a^2}\,+1)),
\endaligned$$
showing (A.4). For (A.5), we furthermore set $u=\sqrt{s^2+1}\,-s$,
noting that $du/ds=(s-\sqrt{s^2+1}\,)/\sqrt{s^2+1}\,$; then
$$\aligned
\int_1^\infty &\tfrac{(\sqrt{r^2+a^2}\,-r)^{j-1}}{
\sqrt{r^2+a^2}\,(\sqrt{r^2+a^2}\,+r)^{j+1}}\, 
dr=a^{-2}\int_{1/a}^\infty \tfrac{(\sqrt{s^2+1}\,-s)^{j-1}}{
\sqrt{s^2+1}\,(\sqrt{s^2+1}\,+s)^{j+1}}\, ds\\
&=a^{-2}\int_{1/a}^\infty \tfrac{(\sqrt{s^2+1}\,-s)^{2j}}{
\sqrt{s^2+1}\,}\, ds=-a^{-2}\int_{s=1/a}^{s=\infty} u^{2j-1} du\\
&=-a^{-2}\bigl[(2j)^{-1}u^{2j}\bigr]_{s=1/a}^{s=\infty } 
=a^{-2}(2j)^{-1}(\sqrt{a^{-2}+1}\,-a^{-1})^{2j},
\endaligned$$ 
since 
$$
\multline
\sqrt{s^2+1}\,-s=s(\sqrt{1+s^{-2}}\,-1)
=s(1+\tfrac12 s^{-2}+O(s^{-4})-1)\\
=\tfrac12 s^{-1}+O(s^{-3})
\text{ for }s\to\infty 
\endmultline
\tag A.6$$
implies $u\to 0$ for $s\to\infty $.
\qed
 
\enddemo

We shall also need:

\proclaim{Lemma A.4} $\int_{-\infty }^\infty
\tfrac{\log(1+s^2)}{1+s^2}\,ds=2\pi \log 2$.
\endproclaim 

\demo{Proof} Write $\log (1+s^2)=\log(1+is)+\log(1-is)$. Here
$\log(1+is)$ extends to a holomorphic function of $s\in \Bbb
C\setminus\{s=it\mid t\ge 1\}$, and $$
\tfrac1{1+s^2}=\tfrac 1{2i}\bigl(\tfrac1{s-i}-\tfrac1{s+i}\bigr).
\tag A.7$$
Since $\log(1+is)/(1+s^2)$ is $O(|s|^{\varepsilon -2})$ for
$|s|\to\infty $, we can calculate$$
\int_{-\infty }^\infty
\tfrac{\log(1+is)}{1+s^2}\,ds=-\int_{\Cal
C_-}\tfrac{\log(1+is)}{1+s^2}\,ds=\tfrac1{2i}\int_{\Cal 
C_-}\tfrac{\log(1+is)}{s+i}\,ds=\pi \log 2,
$$ where we replaced the real axis by a closed curve $\Cal C_-$ in the
lower halfplane going around the 
pole $-i$ in the positive direction, and used (A.7).

The integral with $\log(1-is)$ is turned into an integral over a
curve around the pole $i$ in the upper halfplane, contributing $\pi
\log 2$ in a similar way.\qed
\enddemo

\subhead Acknowledgement \endsubhead 

The author is grateful to  Mirza
Karamehmedovic, Henrik
Schlichtkrull and Elmar Schrohe for helpful discussions of calculations.

\enddocument